\documentclass[12pt, a4paper]{article}

\usepackage[T1]{fontenc}
\usepackage[latin1]{inputenc}
\usepackage[english]{babel}
\usepackage{amsmath}
\usepackage{textcomp}
\usepackage{graphicx}
\usepackage{amsfonts}
\title{KPP reaction-diffusion equations with a non-linear loss inside a cylinder}
\author{Thomas Giletti \thanks{Université Aix-Marseille III, LATP, Faculté des Sciences et Techniques, Avenue Escadrille Normandie-Niemen, F-13397 Marseille Cedex 20, France; thomas.giletti@etu.univ-cezanne.fr}}
\date{2009}

\begin{document}
\maketitle

\begin{abstract}
We consider in this paper a reaction-diffusion system in presence of a flow and under a KPP hypothesis. While the case of a single-equation has been extensively studied since the pioneering Kolmogorov-Petrovski-Piskunov paper, the study of the corresponding system with a Lewis number not equal to 1 is still quite open. Here, we will prove some results about the existence of travelling fronts and generalized travelling fronts solutions of such a system with the presence of a non-linear space-dependent loss term inside the domain. In particular, we will point out the existence of a minimal speed, above which any real value is an admissible speed. We will also give some spreading results for initial conditions decaying exponentially at infinity.
\end{abstract}

\section{Introduction and main results}\label{sec:intro}

There has been a lot of interest in the past years about the effect of flows on the qualitative and quantitative behavior of solutions of reaction-diffusion equations. At first, most of the mathematical analysis only dealt with the flow effect for a single reaction-diffusion equation, studying the existence of travelling fronts \cite{berestycki2,berestycki3,berestycki4,volpert1,xin1,xin2,xin3}, the behavior of the speed of propagation \cite{berestycki0,berestycki5,berestycki6,constantin1,hamel1,heinze1,kiselev1,nolen1}, and flame quenching \cite{constantin2,vladimirova1}. See also \cite{berestycki-hamel,berestycki4,xin3} for reviews of this mathematical area. But recent papers have extended this analysis both in one-dimensional \cite{giovangigli1,roques1,roques2} and multi-dimensional \cite{hamel-quenching,ducrot1,gordon-periodic,hamel-nonadiabatic,hamel-adiabatic} settings to the following system

\begin{equation}\label{eqn:sys1}
\left\{
\begin{array}{l}
T_t + u(y) T_x= \Delta T + f(T)Y, \vspace{3pt} \\
Y_t + u(y) Y_x = \mbox{Le}^{-1} \Delta Y - f(T)Y.\\
\end{array}
\right.
\end{equation}
This problem is posed in a cylinder $\Omega = \mathbb{R}_x \times \omega_y \subset \mathbb{R}^d$ where $\omega$ is a smooth bounded domain of $\mathbb{R}^{d-1}$, with various boundary conditions . We also assume that $u \in C^{0 ,\alpha} (\overline{\omega})$ (for some $\alpha >0$) is the first component of a divergence-free shear flow $(u(y),0)$ with zero average:
\begin{equation}\label{eqn:u}
\int_{\omega} u(y)dy=0 .
\end{equation}
Although this system describes various processes in nature, ranging from chemical and biological contexts (such as predator-prey systems \cite{murray1,murray2}) to combustion and many-particle systems, we will here invoke, to fix the ideas, the "combustion" terminology and refer to the function $T$ as "temperature" and to the function $Y$ as "concentration". The Lewis number Le is then the ratio of the thermal and material diffusivities.

This system is said to be of the KPP-type if $f \in C^1 ([0,+\infty);\mathbb{R})$ and
$$f(0)=0 < f(s) \leq f'(0)s, \ f'(s) \geq 0 \mbox{ for all }s>0 \mbox{ and }f(+\infty)=+\infty .$$
Under this hypothesis, some results have been shown both in the adiabatic case \cite{hamel-adiabatic}, and in the case of heat loss on the boundary \cite{hamel-quenching,hamel-nonadiabatic}, that is with homogeneous Neumann and Robin boundary conditions. See also \cite{gordon-periodic} for some extensions to periodic media with linear heat loss on the boundary. In particular, a sufficient condition has been presented for the existence of travelling fronts, and qualitative properties have been shown to describe such solutions. One of the tools was based on the study of a principal eigenvalue problem in the domain $\omega$, which also permitted to establish criteria for the blow-off and propagation of the solution of the associated Cauchy problem.

\newtheorem{remark}{Remark}
\begin{remark}
\upshape Note that without the KPP hypothesis, the situation is much less clear. For instance, for nonlinearities $f(T)$ of the ignition type (that is, when there exists an ignition temperature $\theta >0$ such that $f(T)=0$ for $T < \theta$ and $f(T)>0$ for $T > \theta $), existence of travelling waves was established only for the Lewis numbers close to 1 in \cite{ducrot1}, or in dimension 1 in \cite{berestycki8,roques1,roques2}.
\end{remark}

In this paper, we will show similar results for the following system, still posed in the cylindric domain $\Omega$:
\begin{equation}\label{eqn:sys}
\left\{
\begin{array}{l}
T_t + u(y) T_x= \Delta T + f(y,T)Y - h(y,T), \vspace{3pt} \\
Y_t + u(y) T_x= \mbox{Le}^{-1} \Delta Y - f(y,T)Y,\\
\end{array}
\right.
\end{equation}
with Neumann boundary conditions
\begin{equation}\label{eqn:neumann}
\frac{\partial T}{\partial n} = \frac{\partial Y}{\partial n} = 0 \mbox{ on } \partial \Omega,
\end{equation}
where $n$ denotes the outward unit normal on $\partial \Omega$. Here, $f\in C^1 (\overline{\omega}\times [0,+\infty) ; \mathbb{R})$ and we assume that there exists $s_0 >0$ such that the set of functions $(f(y,.))_{y\in \overline{\omega}}$ is bounded in $C^{1,\alpha} ([0,s_0) ; \mathbb{R}) $. Moreover, the function $f$ satisfies, by analogy with the KPP case,
$$
\begin{array}{l}
\displaystyle
\displaystyle f(.,0)=0 <  f(.,T) \leq \frac{\partial f}{\partial T}(.,0)T , \ \frac{\partial f}{\partial T} \geq 0 \mbox{ for all } T>0, \ \mbox{and } f(., +\infty )=+\infty ,\\
\end{array}
$$
where the last limit is assumed to be uniform with respect to $y \in \overline{\omega}$.
Furthermore, $h \in C^1 (\overline{\omega}\times [0,+\infty) ; \mathbb{R})$ denotes the heat loss, which takes place in the whole domain, and is such that $(h(y,.))_{y\in \overline{\omega}}$ is bounded in $C^{1,\alpha} ( [0,s_0) ; \mathbb{R})$, along with the conditions
\begin{equation}\label{eqn:condh}
\left\{
\begin{array}{c}
\displaystyle
h(.,0)=0 \leq  \frac{\partial h}{\partial T}(.,0)T \leq h(.,T) \leq K  T < +\infty  \mbox{ for all } T \geq 0 \mbox{ and some } K>0, \\
\displaystyle \int_{\omega} \frac{\partial h}{\partial T}(y,0) dy > 0.
\end{array}
\right.
\end{equation}
For instance, a linear heat loss $h(T) = qT$ where $q>0$ fulfills those hypotheses. The condition on the integral over $\omega$ of $\frac{\partial h}{\partial T} (y, 0)$ means that the heat-loss is non trivially equal to 0 in the domain. It will be used to study the qualitative properties of any solution of (\ref{eqn:sys}). Moreover, the bounds on $h$ are technical hypotheses: $- h(.,T) \leq - \frac{\partial h}{\partial T} (.,0) T$ is similar to the KPP-condition on $f$ and will allow us to use comparisons with the linearized problem, while the boundedness of $\frac{\partial h}{\partial T} $ will allow us to use some standard estimates.

Note also that the space dependence of the heat loss allows us to question whether the solution of (\ref{eqn:sys})-(\ref{eqn:neumann}) converges to a solution of (\ref{eqn:sys1}) with Robin boundary conditions when $h$ converges to a Dirac mass $\delta_{\partial \Omega}$. This will be the subject of a forthcoming paper \cite{giletti2}.

Here, we will follow two main axes. First, we will search for travelling fronts solutions, that is solutions of (\ref{eqn:sys})-(\ref{eqn:neumann}) of the form $T(t,x,y)=\tilde{T} (x-ct,y)$ and $Y(t,x,y)= \tilde{Y} (x-ct, y)$. Thus, we say that $(c,T,Y)$ is a travelling front solution of (\ref{eqn:sys})-(\ref{eqn:neumann}) if in the moving frame $x'=x-ct$ (we drop the primes and the tildes immediately) the functions $T$ and $Y$ satisfy:
\begin{equation}\label{eqn:sysfront}
\left\{
\begin{array}{rcc}
\Delta T + (c-u(y))T_x + f(y,T)Y - h(y,T) & = & 0 \ \\
\mbox{Le}^{-1} \Delta Y + (c-u(y))Y_x - f(y,T)Y & = & 0 \ \\
\end{array} \mbox{ in } \Omega,
\right.
\end{equation}
together with the boundary conditions (\ref{eqn:neumann}) and the following conditions at infinity

\begin{equation}\label{eqn:condinfty}
\left\{
\begin{array}{l}
T(+\infty ,.)=0, \ Y(+\infty ,.)=1 ,\\
T_x (-\infty ,.)=Y_x (-\infty ,.)=0,\\
\end{array}
\right.
\end{equation}
where the limits are uniform with respect to $y\in \overline{\omega}$. The conditions (\ref{eqn:condinfty}) mean that the right-hand side corresponds to the cold region with reactant concentration close to 1, while rather weak conditions are imposed on the left-hand side, that is behind the front. In particular, the values of the temperature and reactant densities are not a priori imposed far behind the front. Furthermore, throughout the paper, the relative concentration $Y$ is assumed to range in $[0,1]$ and is not identically equal to 1. The temperature $T$ is nonnegative and not identically equal to 0.

The other aim of this paper will be to establish criteria for flame blow-off, extinction and propagation. That is, we will consider the solution $(T,Y)$ of the Cauchy problem defined by (\ref{eqn:sys})-(\ref{eqn:neumann}) with an initial profile $(T_0 ,Y_0)$ such that
\begin{equation}\label{eqn:iniprofile}
\begin{array}{l}
\displaystyle 0 \leq T_0 \mbox{, } T_0 \mbox{ is bounded, } 0\leq Y_0 \leq 1 , \\
\displaystyle \exists \lambda >0 \mbox{, } \exists C_1 ,C_2 >0 \mbox{, } C_1 e^{-\lambda x} \leq T_0 (x,y) \leq C_2 e^{-\lambda x} \mbox{ in } \mathbb{R}^+ \times \overline{\omega},\\
\displaystyle \exists \lambda ' >0 \mbox{, } \exists C_3 >0 \mbox{, } 1 - Y_0 (x,y) \leq C_3 e^{-\lambda ' x}  \mbox{ in } \mathbb{R}^+ \times \overline{\omega}.
\end{array}
\end{equation}
We will say that the flame becomes extinct if $\| T(t,.,.) \|_{L^\infty (\Omega )} \rightarrow 0$ as $t \rightarrow +\infty$. The flame is blown-off if there exists a function $\Phi (\xi )$ so that $\Phi (\xi )\rightarrow 0$ as $\xi \rightarrow +\infty$, and $T(t,x,y) \leq \Phi (x+ct )$ with some $c>0$. Lastly, the flame propagates with speed $c>0$ to the right if for any $c' >c$, $T(t, x+c' t,y) \rightarrow 0$ as $t \rightarrow +\infty$ while for the speed $c$ itself, one can find $x_0 \in \mathbb{R}$ and $\alpha (x_0, y) > 0$ such that $T(t,x_0 + ct ,y) \geq \alpha (x_0, y)$ for all $t \geq 1$ and $y \in \overline{\omega}$.

Before we state the main results of this paper, we introduce the following principal eigenvalue problem depending on a parameter $\lambda \in \mathbb{R}$:

\begin{equation}\label{eqn:principaleigenvalueh}
\left\{
\begin{array}{rcll}
\displaystyle -\Delta_y \phi_\lambda - \lambda u(y) \phi_\lambda + ( \frac{\partial h}{\partial T} (y,0) -\frac{\partial f}{\partial T} (y,0))  \phi_\lambda & = & \mu_{h,f} (\lambda ) \phi_\lambda & \mbox{ in } \omega ,\\
\displaystyle \frac{\partial \phi_\lambda }{\partial n} & = & 0 & \mbox{ on } \partial \omega .\\
\end{array}
\right.
\end{equation}
That is, $\mu_{h,f} (\lambda )$ is the unique eigenvalue of (\ref{eqn:principaleigenvalueh}) that corresponds to a positive eigenfunction $\phi_\lambda (y)$, and can be defined for any functions $f$, $h \in C^1 (\overline{\omega} \times [0, +\infty ); \mathbb{R})$. Let us first show some properties of the function $\mu_{h,f}$. The eigenfunction $\phi_\lambda$ can be normalized so that
$$\int_\omega \phi_\lambda ^2 (y) dy =1.$$
With this normalization, one gets
\begin{equation}
\mu_{h,f} (\lambda ) = \int_\omega | \nabla \phi_\lambda (y) |^2 dy - \lambda \int_\omega u(y) \phi_\lambda ^2 (y) dy + \int_\omega ( \frac{\partial h}{\partial T} (y,0) -\frac{\partial f}{\partial T} (y,0)) \phi_\lambda ^2 (y) dy.
\end{equation}
In particular, under the hypotheses made on $h$, we have that
$$\int_\omega \frac{\partial h}{\partial T} (y,0) \phi_\lambda ^2 (y) dy \geq \min_{\overline{\omega}} \phi_\lambda ^2 \times \int_\omega \frac{\partial h}{\partial T} (y,0)  dy >0,$$
which implies that $\mu_{h,0} (0) >0$. Furthermore, by the variational principle, we have that
\begin{eqnarray*} \displaystyle
\mu_{h,f} (\lambda ) = \min_{\psi \in H^1 (\omega ), \| \psi  \|_2=1} \displaystyle & \displaystyle \left(
\int_\omega | \nabla \psi (y) |^2 dy - \lambda \int_\omega u(y) \psi ^2 (y) dy \right. \\
& \displaystyle \left. \ \ \ \  \ + \int_\omega ( \frac{\partial h}{\partial T} (y,0) -\frac{\partial f}{\partial T} (y,0)) \psi ^2 (y) dy \right) ,
\end{eqnarray*}
where $\| . \|_2$ denotes the $L^2 (\omega )$ norm. This implies that $\mu_{h,f} (\lambda)$ is concave as an infimum of a family of affine functions. We now give one last property of $\mu_{h,f}$, which will allow us to discuss the conditions of our theorems later in this paper.

Remember first that $\lambda \mapsto \mu_{h,f} (\lambda )$ and $\lambda \mapsto \phi_\lambda$ are analytic functions of $\lambda$. When differentiating (\ref{eqn:principaleigenvalueh}) with respect to $\lambda$, we obtain
\begin{equation}\label{eqn:principaleigenvaluehprime}
-\Delta \phi_\lambda ' - \lambda u(y) \phi_\lambda ' - u(y) \phi_\lambda + ( \frac{\partial h}{\partial T} (y,0) -\frac{\partial f}{\partial T} (y,0)) \phi_\lambda ' = \mu_{h,f} ' (\lambda ) \phi_\lambda + \mu_{h,f} (\lambda ) \phi_\lambda ' ,
\end{equation}
where the prime denotes derivative with respect to $\lambda$. By multiplying (\ref{eqn:principaleigenvalueh}) by $\phi_\lambda '$, (\ref{eqn:principaleigenvaluehprime}) by $\phi_\lambda$, substracting one equation from the other and using the $L^2$-normalization of $\phi_\lambda$, we obtain that
\begin{equation}\label{eqn:muprime}
\mu_{h,f} ' (\lambda) = - \int_\omega u(y) \phi_\lambda^2 (y) dy.
\end{equation}
We also introduce the following principal eigenvalue problem, also depending on a parameter $\lambda \in \mathbb{R}$:
\begin{equation}\label{eqn:principaleigenvalue}
\left\{
\begin{array}{rcll}
\displaystyle  -\Delta_y \psi_\lambda - \lambda u(y) \psi_\lambda & = & \nu (\lambda ) \psi_\lambda & \mbox{ in } \omega ,\\
\displaystyle \frac{\partial \psi_\lambda }{\partial n} & = & 0 & \mbox{ on } \partial \omega .\\
\end{array}
\right.
\end{equation}
That is, $\nu (\lambda )$ is the unique eigenvalue of (\ref{eqn:principaleigenvalue}) that corresponds to a positive eigenfunction $\psi_\lambda (y)$. In fact, this is the same principal eigenvalue problem, with $h=f=0$ (the purpose of its introduction is only to simplify some of our notations). In particular, we can obtain as above that $\nu (\lambda)$ is concave. Furthermore, (\ref{eqn:u}) and (\ref{eqn:muprime}) with $h=f=0$, together with the fact that any positive constant is an eigenfunction of (\ref{eqn:principaleigenvalue}) with $\lambda =0$, imply that $\nu (0) =\nu ' (0) =0$. This in turn implies that $\nu (\lambda)$ is nonpositive for all $\lambda \in \mathbb{R}$.
\begin{remark}\label{rem:hy}
\upshape When $f(y,T)=f(T)$ and $h(y,T) = h(T)$ do not depend on $y$, it is immediate to see that $\mu_{h,f} (\lambda ) = \nu (\lambda) + h'(0) - f'(0)$. In particular, $\mu_{h,f} (\lambda ) \leq h'(0) - f'(0)= \mu_{h,f} (0)$ for all $\lambda \in \mathbb{R}$ in that case.
\end{remark}

In Section \ref{sec:qualit}, we will show the following qualitative properties of any travelling front solution:
\newtheorem{theorem}{Theorem}
\begin{theorem}\label{th:thqualit}
Let $(c,T,Y)$ be a solution of $(\ref{eqn:sysfront})$-$(\ref{eqn:condinfty})$ and $(\ref{eqn:neumann})$ such that $0<T$ and $0<Y<1$. Then $T$ is bounded, $T(-\infty ,.)=0$, $Y(-\infty ,.)=Y_{\infty} \in (0,1)$, $\mu_{h,f} (0) <0$, $c>0$ and $c\geq c^*$,
where $c^*$ is then defined by:\begin{equation}\label{eqn:c*}
c^* = \min \{c \in \mathbb{R}, \ \exists \lambda >0 , \ \mu_{h,f} (\lambda) = \lambda^2 - c \lambda \} = \min_{\lambda >0} \frac { k(\lambda)}{\lambda},
\end{equation}
and $$k (\lambda ) = \lambda^2 -\mu_{h,f} (\lambda) .$$
\end{theorem}
Note that since $k$ is strictly convex and under the hypothesis $\mu_{h,f} (0) <0$, it is straightforward to check that the equation $k(\lambda ) =c\lambda$ has one positive solution $\lambda^*$ for $c=c^*$, and two positive solutions $\lambda_1$, $\lambda_2$ for $c>c^*$, with $\lambda_1 < \lambda^* <\lambda_2$.

We will then show the existence of bounded travelling fronts solutions:
\begin{theorem}\label{th:thexist}
$(a)$ Assume that $\mu_{h,f} (0) < 0$. For any $c > \max (0, c^*)$, there exists a solution $(T,Y)$ of $(\ref{eqn:sysfront})$-$(\ref{eqn:condinfty})$ and $(\ref{eqn:neumann})$ such that $T$ is bounded, $T(-\infty ,.)=0$, $T>0$, $0<Y<1$ and $Y(-\infty ,.) = Y_{\infty } \in (0,1)$.\\
$(b)$ Assume that $\sup_{\lambda \in \mathbb{R}} \displaystyle (\mu_{h,f} (\lambda ) - \lambda^2 )< 0$. Then $c^* >0$ and there exists a solution $(T,Y)$ of $(\ref{eqn:sysfront})$-$(\ref{eqn:condinfty})$ and $(\ref{eqn:neumann})$ with minimal speed $c=c^*$, and such that $T$ is bounded, $T(-\infty ,.)=0$, $T>0$, $0<Y<1$ and $Y(-\infty ,.) = Y_{\infty } \in (0,1)$.
\end{theorem}

\begin{remark}\label{rem:minspeed1}
\upshape It immediately follows from Remark \ref{rem:hy} that when $h$ and $f$ are independent of $y \in \overline{\omega}$ and under the hypothesis $\mu_{h,f} (0) = h'(0)- f'(0) <0$, both parts of Theorem~\ref{th:thexist} are verified. That is, there exists a non trivial travelling front solution for any speed $c \geq c^* >0$.
\end{remark}

Section \ref{sec:exist} will be dedicated to the proof of part (a). Part (b) will be treated in Section~\ref{sec:existc*}. Lastly, Section \ref{sec:cauchy} will deal with the Cauchy problem, with a proof of the following result.

\begin{theorem}\label{th:cauchy}
Let $(T,Y)$ be a solution of $(\ref{eqn:sys})$-$(\ref{eqn:neumann})$ with an initial profile $(T_0 ,Y_0)$ verifying $(\ref{eqn:iniprofile})$. Let $\lambda$ be the decay rate of $T_0$ as in $(\ref{eqn:iniprofile})$.\\
$(a)$ Extinction. If $\mu_{h,f} (0) > 0$, then $T(t,x,y) \leq C e^{-\gamma t}$ for all $t\geq 0$ and $(x,y) \in \overline{\Omega}$ where $ \gamma = \mu_{h,f} (0) >0$ and $C$ is a positive constant.\\
$(b)$ Blow-off. Let us assume that there exists \ $0 < \eta \leq \lambda$ such that $\mu_{h,f} (\eta) - \eta^2 > 0 $. Then $T(t,x,y) \leq C e^{-\eta (x+ \delta t)}$ for all $t \geq 0$ and $(x,y) \in \overline{\Omega}$ with $C$, $ \delta >0$.\\
$(c)$ Propagation. Let us assume that $\mu_{h,f} (0) < 0$, $\mu_{h,f} (\lambda) -\lambda^2 < 0$ and $\lambda < \lambda^*$. Then $c := k(\lambda ) / \lambda > \max(0,c^*)$ and the solution propagates with speed $c$.
\end{theorem}
Parts (a) of Theorem \ref{th:thexist} and \ref{th:cauchy} reflects the fact that $\mu_{h,f} (0) < 0$ is a sufficient condition for the existence of a travelling front solution, and shows that it is also almost a necessary condition for the propagation of the flame (the case $\mu_{h,f} (0) = 0$ is still open). It is also important to note that Part (c) of Theorem~\ref{th:cauchy} underlines the link between the speed of propagation and the decay rate of temperature on the right, which will in fact be used several times throughout the paper.

Lastly, let us discuss the completeness of this Theorem. If $\mu_{h,f} (0)> 0$, we can apply part (a). If $\mu_{h,f} (0) <0$, we first consider the case $c^* < 0$. Let then $\lambda_1 < \lambda_2$ the solutions of $\mu_{h,f} (s) - s^2 =0$. We have $\lambda_1 < \lambda^* < \lambda_2$. Furthermore, $\mu_{h,f} (s) - s^2 $ is negative for $s \in (0, \lambda_1)$ and positive for $s \in (\lambda_1 ,\lambda_2)$. Thus, part (b) can apply for $\lambda > \lambda_1$, and part (c) for $\lambda < \lambda_1$. In the case $c^* \geq 0$, we have that $\mu_{h,f} (s) - s^2 $ is negative for $s \in (0, \lambda^*)$ and nonpositive everywhere. Thus, part (c) apply for $\lambda < \lambda^*$, but the problem is still open for $\lambda \geq \lambda^*$.
In the latter case, we may at least say that the solution can't propagate with speed $c>c^*$, by placing ourselves in a moving frame with speed $c$, and then using part(b) of our Theorem. This argument, along with the well-known fact of the propagation with minimal speed in the single-equation case for a heaviside initial condition \cite{berestycki1}, may allow us to conjecture that the solution propagates with speed $c^*$ for $\lambda>\lambda^*$. Nevertheless, at this time, no significant result has been made in this direction to our knowledge in the system case.

\section{Qualitative properties of travelling fronts} \label{sec:qualit}
This section is devoted to the proof of Theorem \ref{th:thqualit}. Let $(c,T,Y)$ be a solution of (\ref{eqn:sysfront})-(\ref{eqn:condinfty}) and (\ref{eqn:neumann}) such that $0<T$ and $0<Y<1$.
\subsection{Boundedness of temperature} \label{sec:Tbounded}

We first prove that $Y$ converges to a constant as $x \rightarrow -\infty$. To this end, we integrate equation (\ref{eqn:sysfront}) satisfied by $Y$ over the domain $(-N,N) \times \omega$ with $N>0$. We obtain:
\begin{equation}\label{eqn:eqint1}
\begin{array}{c}
\displaystyle \int_\omega [ \mbox{Le}^{-1} (Y_x (N,y) - Y_x (-N,y)) + (c-u(y))(Y(N,y)-Y(-N,y))]dy\\
\\
=\displaystyle \int_{(-N,N) \times \omega} f(y,T(x,y))Y(x,y) dx dy.
\end{array}
\end{equation}
Recall that $Y_x (x,y) \rightarrow 0$ as $x\rightarrow - \infty$ from (\ref{eqn:condinfty}). Besides, as $T$ is bounded for $x>0$ and $Y$ converges to $1$ as $x \rightarrow +\infty$, it follows from standard elliptic estimates that $Y_x (x,y) \rightarrow 0$ as $x \rightarrow +\infty$. Since $Y$ is bounded, we finally have that the left-hand side is bounded independently of $N$. Therefore, as $f(y,T)Y$ is a positive function, we conclude that the positive integral
$$\int_{\Omega} f(y,T(x,y))Y(x,y) dxdy$$ converges. Furthermore, by multiplying the equation (\ref{eqn:sysfront}) satisfied by $Y$ by $Y$ itself, and integrate over the domain $(-N,N) \times \omega$ with $N>0$, we obtain:
\begin{equation*}
\begin{array}{c} \displaystyle \int_\omega [ \mbox{Le}^{-1} (Y_x (N,y) Y(N,y) - Y_x (-N,y)Y(-N,y)) + \frac{1}{2} (c-u(y))(Y^2 (N,y) -Y^2 (-N,y))]dy \\
\\
=\displaystyle \int_{(-N,N) \times \omega} [f(y,T(x,y))Y^2 (x,y) + \mbox{Le}^{-1} |\nabla Y | ^2 ]dx dy.
\end{array}
\end{equation*}
The left-hand side is again bounded independently of $N$, and so is the integral
$$0 \leq \int_{(-N,N) \times \omega} f(y,T(x,y))Y^2 (x,y) dx dy < +\infty.$$
We conclude that the integral
$$\int_{\Omega} |\nabla Y(x,y)|^2 dxdy$$
converges. Still, since the function $T$ is not known to be a priori bounded, we can't use $W_{loc}^{2,p}$ estimates to prove the convergence of $Y$ to a constant as $x \rightarrow +\infty$. To overcome this difficulty, we fix $a \in \mathbb{R}$ and let $(x_k)_{k \in \mathbb{N}}$ be any sequence converging to $-\infty$ as $k\rightarrow +\infty$. We introduce the translate $Y_k (x,y)=Y(x_k +a+x,y)$. We have
$$\int_{(a+x_k ,a+1+x_k) \times \omega} | \nabla Y |^2 dxdy = \int_{(0,1)\times \omega} | \nabla Y_k |^2 dxdy \ \rightarrow \ 0 \mbox{ as } k\rightarrow +\infty .$$
Hence, up to extraction of a subsequence, $(Y_k)_{k \in N}$ converges in $H^1 ((0,1)\times \omega)$ to a constant $Y_\infty^a \in [0,1]$. We then use (\ref{eqn:eqint1}) with $N=-x_k -a-\xi$ for $\xi \in (0,1)$ ($k$ is chosen large enough so that $-x_k - a-1 >0$) and integrate over $\xi \in (0,1)$. We obtain
$$\int_{(0,1)\times \omega} \mbox{Le}^{-1}(Y_x (-x_k -a-\xi ,y)-Y_x (x_k+a+\xi ,y))d\xi dy$$
$$+\int_{(0,1)\times \omega} (c-u(y))(Y(-x_k -a-\xi ,y)-Y(x_k +a+\xi ,y))dy$$
$$=\int_0^1
\left(
\begin{array}{l}
\displaystyle \int_{(x_k +a+\xi ,-x_k -a-\xi )\times \omega} f(y,T(x,y))Y(x,y)dxdy
\end{array}
\right)
d\xi .$$
The first term of the left side converges to $0$ as $k \rightarrow +\infty$ (recall that $Y_x (+\infty ,.)=Y_x (-\infty ,.) =0$ uniformly in $y \in \overline{\omega}$). The second term of the left side converges to
$$\int_\omega (c- u(y))(1-Y_\infty^a ) \ = \ (1-Y_\infty^a )c| \omega | .$$
We used here the fact that $u(y)$ has mean zero (\ref{eqn:u}). Lastly, by the dominated convergence theorem, the right-hand side converges to
$$\int_{\Omega} f(y,T(x,y))Y(x,y) dx dy .$$
Therefore,
\begin{equation}\label{eqn:Yinftya}
(1-Y_\infty^a )c | \omega | = \int_{\Omega} f(y,T(x,y))Y(x,y) dxdy >0 ,
\end{equation}
and as a consequence, $Y_\infty^a < 1$ does not depend on $a$ nor on the sequence $(x_k)_{k \in \mathbb{N}} $. It also already implies that $c>0$. We conclude that there exists a constant $Y_\infty \in [0,1)$ such that $Y(x_k +. ,.)$ converges to $Y_\infty$ in $H_{loc}^1 (\Omega )$ as $k\rightarrow +\infty$ for any sequence $(x_k)_{k \in \mathbb{N}} \rightarrow -\infty$.
\\

Let us now prove that $T$ is globally bounded. Assume by contradiction that $T$ is unbounded. Since $T(+\infty ,.)=0$, there has to exist a sequence $(x_k ,y_k )_{k\in \mathbb{N}}$ in $\mathbb{R} \times \omega$ such that $x_k \rightarrow -\infty$ and
\begin{equation}\label{eqn:lim1}
T(x_k ,y_k) \rightarrow +\infty
\end{equation}
as $k \rightarrow +\infty$. Since the functions $Y$, $f(y,T)/T$ and $h(y,T))/T$ are bounded in $\Omega$, it follows from standard elliptic estimates and the Harnack inequality up to the boundary that $|\nabla T|/T$ is also bounded in $\Omega$. Thus, we also have
\begin{equation}\label{eqn:lim2}
\min_{(x,y)\in [x_k -1 , x_k +1] \times \overline{\omega}} T(x,y) \rightarrow +\infty
\end{equation}
as $k\rightarrow +\infty$. Hence, we also have
$$\min_{(x,y)\in [x_k -1 , x_k +1] \times \overline{\omega}} f(y,T(x,y)) \rightarrow +\infty$$
as $k\rightarrow +\infty$, because $f(y,+\infty ) = +\infty$ uniformly in $y \in \overline{\omega}$. But
$$\int_{\Omega} f(y,T(x,y))Y(x,y)dxdy \geq \int_{(x_k -1 ,x_k +1) \times \omega} f(y,T(x,y))Y(x,y)dxdy$$
$$\geq \min_{(x,y) \in [x_k -1 , x_k +1] \times \overline{\omega}} f(y,T(x,y)) \ \times \ \int_{(x_k -1 ,x_k +1) \times \omega} Y(x,y)dxdy,$$
and $$ \int_{(x_k -1 ,x_k +1) \times \omega} Y(x,y)dxdy \rightarrow 2 | \omega | Y_\infty$$
as $k \rightarrow +\infty$. We conclude that if $T$ is unbounded, then $Y_\infty =0$. Let now introduce the functions
$$T_k (x,y) =\frac{T(x_k +x ,y)}{T(x_k ,y_k )},$$
which are locally bounded, as follows from the fact that $| \nabla T | /T \in L^\infty (\overline{\Omega})$. These functions satisfy
$$
\left\{
\begin{array}{rcll}
\displaystyle \Delta T_k + (c-u(y))T_{k,x} + (g_{1,k} -g_{2,k})  T_k & = & 0 & \mbox{in } \Omega ,\\
\displaystyle \frac{\partial T_k}{\partial n} & = & 0 & \mbox{on } \partial \Omega ,\\
\end{array}
\right.
$$
where
$$\displaystyle g_{1,k} (x,y) = \frac{f(y,T(x_k +x ,y))}{T(x_k +x ,y)} Y (x_k +x ,y) ,$$
$$\displaystyle g_{2,k} (x,y) = \frac{h(y,T(x_k +x ,y))}{T(x_k +x ,y)} \ .$$
First, we have
$$\displaystyle 0 \leq g_{1,k} (x,y) \leq \frac{\partial f}{\partial T} (y,0) \ Y (x_k +x ,y)  \leq \max_{\overline{\omega}} \frac{\partial f}{\partial T} (.,0),$$
and $g_{1,k} \rightarrow 0$ in $L_{loc}^2 (\overline{\Omega})$ because $Y_\infty =0$. On the other hand,
$$0 \leq \frac{\partial h}{\partial T} (y,0)  \leq g_{2,k} (x,y) \leq  K ,$$
where $K$ comes from the bounds on $h$ from (\ref{eqn:condh}).
Thus, as the sequence $(g_{2,k})_{k \in \mathbb{N}}$ is bounded in $L^\infty (\overline{\Omega})$, it converges weakly in $L^{1,*} (\overline{\Omega})$ up to extraction of some subsequence to a function $g$  in $L^\infty (\overline{\Omega})$.

Lastly, since the functions $g_{1,k}$ and $g_{2,k}$ are uniformly bounded in $L^\infty (\Omega)$, the functions $T_k$ are then bounded in $W_{loc}^{2,p} (\overline{\Omega})$ for all $1\leq p < +\infty$. Up to extraction of a subsequence, they then converge weakly in $W_{loc}^{2,p} (\overline{\Omega})$ for all $1\leq p < +\infty$ and then in $C_{loc}^{1,\beta} (\overline{\Omega})$ for all $0 \leq \beta <1$, to a nonnegative solution $T_\infty $ of
$$\left\{
\begin{array}{rcll}
\displaystyle \Delta T_\infty + (c-u(y)) T_{\infty ,x} - g T_\infty & = & 0 & \mbox{in } \Omega ,\\
\displaystyle \frac{\partial T_\infty}{\partial n} & = & 0 & \mbox{on } \partial \Omega .\\
\end{array}
\right.
$$
The elliptic regularity theory implies that the function $T_\infty$ is actually of the class $C_{loc}^{2,\alpha} (\overline{\Omega})$ (remember that $u \in C^{0, \alpha} (\overline{\omega})$). It follows from the boundary condition $T_x \rightarrow 0$ as $x \rightarrow -\infty$ in (\ref{eqn:condinfty}) and from (\ref{eqn:lim1}) that $T_{k,x} (x,y) \rightarrow 0$ locally uniformly as $k\rightarrow +\infty$, whence $T_\infty$ is a function of $y$ only. It is then a solution of
\begin{equation}\label{eqn:sysTinfty1}
\left\{
\begin{array}{rcll}
\displaystyle \Delta T_\infty - g T_\infty & = & 0 & \mbox{in } \Omega ,\\
\displaystyle \frac{\partial T_\infty}{\partial n} & = & 0 & \mbox{on } \partial \Omega .
\end{array}
\right.
\end{equation}
Furthermore, $T_k (0,y_k)=1$ and one can assume, up to extraction of another subsequence, that the sequence $y_k$ converges to $y_\infty \in \overline{\omega}$ as $k\rightarrow +\infty$. Therefore, $T_\infty (0,y_\infty ) =1$ and the strong maximum principle and the Hopf lemma imply that $T_\infty$ is positive in $\overline{\Omega}$. Here, recall that $g$ is the limit in $L^{1,*} (\overline{\Omega})$ of the sequence $(g_{2,k})_{k\in \mathbb{N}}$ where $g_{2,k} (x,y) \geq \frac{\partial h}{\partial T}(y,0)$ for all $(x,y) \in \Omega$. Then, for any $N>0$, we have that
\begin{equation}\label{eqn:contradTbounded}
\begin{array}{rcl}
\displaystyle \int_{(-N ,N) \times \omega } g (x,y) T_\infty (x,y) dx dy & \geq & \displaystyle \int_{(-N ,N) \times \omega } \frac{\partial h}{\partial T}(y,0) T_\infty (x,y) dxdy \\
& \geq & \displaystyle 2N \int_{\omega } \frac{\partial h}{\partial T}(y,0)dy \times \min_{(-N ,N) \times  \overline{\omega}} T_\infty \\
& > & 0.
\end{array}
\end{equation}
However, since $T_{\infty ,x} =0$ and because of the Neumann boundary condition on $\partial \Omega$, integrating (\ref{eqn:sysTinfty1}) over $(-N ,N) \times \omega$ leads to
$$ \int_{ (-N, N) \times \omega} g(x,y) T_\infty (x,y) dxdy =0.$$
This enters in contradiction with (\ref{eqn:contradTbounded}). We conclude that $T$ belongs to $L^\infty (\Omega)$.

\subsection{The left limit for temperature}\label{sec:leftlimitT}

We now show that $T \rightarrow 0$ as $x \rightarrow -\infty$. Recall that the integral
$$\int_{\Omega} f(y,T(x,y))Y(x,y) dxdy$$ converges. We now integrate the equation (\ref{eqn:sysfront}) satisfied by $T$ over the domain $(-N ,N) \times \omega$ with $N>0$. We obtain
\begin{equation}\label{eqn:eqintTb}
\begin{array}{c}
\displaystyle \int_\omega [(T_x (N,y) - T_x (-N,y)) + (c-u(y))(T(N,y) - T(-N,y))]dy \\
\displaystyle =\int_{(-N,N) \times \omega} h(y,T(x,y))dxdy -\int_{(-N,N) \times \omega} f(y,T(x,y))Y(x,y) dxdy .
\end{array}
\end{equation}
It follows from (\ref{eqn:condinfty}) and standard elliptic estimates that $T_x (\pm \infty ,.)=0$. Since $T \in L^\infty (\Omega)$, we deduce that the left-hand side is bounded independently of $N$, whence the integral
$$\int_{\Omega} h(y,T(x,y)) dxdy$$
converges. Besides, since $T$ is bounded and $h(y,T) \geq \frac{\partial h }{\partial T} (y,0) T \geq 0$, we also have that:
\begin{equation}\label{eqn:intconv1} \int_{\Omega} \frac{\partial h }{\partial T} (y,0) T(x,y) dxdy \ < +\infty ,
\end{equation}
$$\int_{\Omega} h(y,T(x,y))T(x,y) dxdy \ < \ +\infty ,$$
$$\int_{\Omega} f(y,T(x,y))Y(x,y)T(x,y) dxdy \ < \ +\infty .$$
We now multiply the equation (\ref{eqn:sysfront}) satisfied by $T$ by $T$ itself and integrate over the domain $(-N ,N) \times \omega$ with $N>0$. We obtain
$$\int_\omega [(T_x (N,y)T(N,y) - T_x (-N,y)T(-N,y)) + (c-u(y))(T^2 (N,y) - T^2(-N,y))]dy$$
$$=\int_{(-N,N) \times \omega} h(y,T)Tdxdy -\int_{(-N,N) \times \omega} f(y,T)YT dxdy + \int_{(-N,N) \times \omega} | \nabla T |^2 dxdy .$$
As before, the left-hand side is bounded independently of $N$ and we saw that the first two integrals of the right-hand side converge as $N\rightarrow +\infty$, whence
\begin{equation}\label{eqn:intconv2}
\int_{\Omega} | \nabla T |^2 dxdy \ < \ +\infty .
\end{equation}
Let now $(x_k)_{k \in \mathbb{N}}$ be any sequence such that $x_k \rightarrow -\infty$ as $ k\rightarrow +\infty$. We define the functions $T_k (x,y) = T(x+x_k ,y)$ for each $k \in \mathbb{N}$. It follows from standard elliptic estimates that this sequence is bounded in $W_{loc}^{2,p} (\overline{\Omega} )$ for all $1 \leq p < +\infty$. Therefore, up to extraction of a subsequence, it converges in $C_{loc}^{1} (\overline{\Omega} )$ to a function $T_\infty$. Because of (\ref{eqn:intconv2}), we know that $T_\infty$ is a constant. Furthermore, it follows from (\ref{eqn:intconv1}) that
$$\int_{(x_k -1 , x_k +1 ) \times \omega} \frac{\partial h }{\partial T} (y,0) T(x,y) dxdy  \rightarrow 2 T_\infty \int_{\omega} \frac{\partial h }{\partial T} (y,0) dy =0$$
as $k \rightarrow + \infty$, whence $T_\infty =0$ does not depend on the choice of the sequence $(x_k)_{k \in \mathbb{N}}$. We conclude that $T(x,y) \rightarrow 0$ when $x \rightarrow -\infty$ locally uniformly in $y \in \overline{\omega}$.

\subsection{Proof of the inequality : $\mu_{h,f} (0) <0$}

Assume by contradiction that $\mu_{h,f} (0) \geq 0$. Let $\phi = \phi_0 $ a positive solution of (\ref{eqn:principaleigenvalueh}) with $\lambda =0$. The function $\phi$ satisfies
\begin{equation}\label{eqn:fmucomp1}
\left\{
\begin{array}{rl}
\displaystyle \Delta \phi +  ( \frac{\partial f}{\partial T} (y,0)- \frac{\partial h}{\partial T} (y,0)) \phi \leq \Delta \phi + \mu_{h,f} (0) \phi + ( \frac{\partial f}{\partial T} (y,0)- \frac{\partial h}{\partial T} (y,0)) \phi=0 & \mbox{in } \omega ,\\
\displaystyle \frac{\partial \phi}{\partial n}  =  0 & \mbox{on } \partial \omega .\\
\end{array}
\right.
\end{equation}
Since $T$ is globally bounded and $\phi$ positive on $\overline{\omega}$, there exists $\gamma >0$ such that $T(x,y) \leq~\gamma \phi (y)$ in $\overline{\Omega}$. Since $T>0$ and $T(\pm \infty ,.)=0$, there exists then $\gamma^* >0$ such that $T (x,y) \leq \gamma^* \phi (y)$ in $\overline{\Omega}$ with equality somewhere. But since $T >0$ and $Y<1$, the function $T$ satisfies
$$\Delta T + (c-u(y)) T_x + \frac{\partial f}{\partial T} (y,0)T - \frac{\partial h}{\partial T} (y,0) T \geq
\Delta T + (c-u(y)) T_x + f(T)Y - h(y,T) =0$$
in $\Omega$ with the Neumann boundary conditions on $\partial \Omega$. Let now $z (x,y)= T (x,y) - \gamma^* \phi (y)$. $z$ is nonpositive in $\overline{\Omega}$ and vanishes somewhere. Besides, $z$ satisfies
$$\Delta z + (c-u(y)) z_x + (\frac{\partial f}{\partial T} (y,0)- \frac{\partial h}{\partial T} (y,0)) z \geq 0$$
in $\Omega$, together with the Neumann boundary conditions on $\partial \Omega$. It then follows from the strong maximum principle and the Hopf lemma that $z=0$, whence $T (x,y) = \gamma^* \phi (y)$ in $\Omega$. This is impossible since $\gamma^* >0$ and $T(+\infty ,.)=0$. We then conclude that $\mu_{h,f} (0) < 0$.

\subsection{A lower bound for the front speed}

Recall that we already saw in Section \ref{sec:Tbounded} that $c >0$. We now prove that $c \geq c^*$ where $c^*$ defined in Section \ref{sec:intro}.

Recall that from Harnack inequality, $|\nabla T|/T$ is also bounded in $\Omega$. Let $\Lambda$ be defined by:
$$\Lambda = - \liminf_{x \rightarrow +\infty}
\left(
\begin{array}{l}
\displaystyle  \min_{y \in \overline{\omega}} \frac{T_x (x,y)}{T(x,y)}\\
\end{array}
\right)
.
$$
Since $T>0$, and $T(+\infty ,.)=0$, we have $\Lambda \geq 0$. Now let $(x_n ,y_n)_{n \in \mathbb{N}} $ be a sequence of points in $\overline{\Omega}$, such that $x_n \rightarrow +\infty$ and
$$\frac{T_x (x_n , y_n)}{T(x_n ,y_n)} \rightarrow \Lambda \mbox{ as } n \rightarrow +\infty .$$
Up to extraction of a subsequence, one can assume that $y_n \rightarrow y_\infty \in \overline{\omega}$ as $n \rightarrow +\infty$. Next, define the normalized and shifted temperature for all $n \in \mathbb{N}$ and $(x,y) \in \overline{\Omega}$:
$$T_n (x,y) = \frac{T(x+x_n ,y)}{T(x_n ,y_n)}.$$
Since $|\nabla T|/T$ is bounded in $\Omega$, the sequence of functions $T_n$ is bounded in $L_{loc}^\infty (\overline{\Omega})$. We also have for each $n \in \mathbb{N}$, $T_n$ is a solution of the following problem:
$$
\left\{
\begin{array}{rcll}
\displaystyle \Delta T_n + (c-u(y))T_{n,x} + \frac{f(y,T(x_n,y_n)T_n)Y_n - h(y, T(x_n ,y_n)T_n)}{T(x_n,y_n)} & = & 0 & \mbox{in } \Omega ,\\
\displaystyle \frac{\partial T_n}{\partial n} & = & 0 & \mbox{on } \partial \Omega ,\\
\end{array}
\right.
$$
where $$Y_n (x,y)=Y(x+x_n ,y)$$
is the shifted concentration.

Recall that $T(x+x_n,y) \rightarrow 0$ and $Y(x+x_n ,y) \rightarrow 1$ locally uniformly in $(x,y) \in \overline{\Omega}$ as $n\rightarrow +\infty$ because of (\ref{eqn:condinfty}). It then follows from standard elliptic estimates that, up to extraction of a subsequence, the sequence $T_n$ converges weakly in $W_{loc}^{2,p} (\overline{\Omega})$ for all $1 \leq p< +\infty$ and strongly in $C_{loc}^{1,\beta } (\overline{\Omega})$ for all $0 \leq \beta <1$, to a function $T_\infty$ which satisfies
$$
\left\{
\begin{array}{rcll}
\displaystyle \Delta T_\infty + (c-u(y))T_{\infty,x} + (\frac{\partial f}{\partial T} (y,0) -  \frac{\partial h}{\partial T}(y,0)) T_\infty & = & 0 & \mbox{in } \Omega , \\
\displaystyle \frac{\partial T_\infty}{\partial n} & = & 0 & \mbox{on } \partial \Omega .\\
\end{array}
\right.
$$
Since $T_n (x,y) \geq 0$ and $T_n (0,y_n) =1$, we also have $T_\infty \geq 0$ in $\overline{\Omega}$ and $T_\infty (0,y_\infty) =1$, whence $T_\infty >0$ in $\overline{\Omega}$, as follows from the Hopf lemma and the strong maximum principle. We can then define $z=T_{\infty ,x}/T_\infty$, which satisfies $$z \geq -\Lambda \mbox{ in } \overline{\Omega},$$
and $z(0,y_\infty ) = -\Lambda$ owing to the definition of $\Lambda$ and the choice of the sequence $(x_n ,y_n )$. Moreover, the function $z$ satisfies the following elliptic equation:
$$\left\{
\begin{array}{rcll}
\displaystyle \Delta z + 2 \frac{\nabla T_\infty}{T_\infty}.\nabla z + (c-u(y))z_x & = & 0 & \mbox{in } \Omega ,\\
\displaystyle \frac{\partial z}{\partial n} & = & 0 & \mbox{on } \partial \Omega .\\
\end{array}
\right.
$$
It is then implied by Hopf lemma and the strong maximum principle that $z(x,y)=-\Lambda$ in $\overline{\Omega}$, that is $T_\infty (x,y) = e^{-\Lambda x} \phi (y)$ in $\overline{\Omega}$ where $\phi (y)$ is a positive function and satisfies
$$\left\{
\begin{array}{rcll}
\displaystyle - \Delta_y \phi - \Lambda u(y) \phi + ( \frac{\partial h}{\partial T} (y,0) - \frac{\partial f}{\partial T} (y,0)) \phi & = & (\Lambda^2 -\Lambda c ) \phi & \mbox{in } \omega ,\\
\displaystyle \frac{\partial \phi}{\partial n} & = & 0 & \mbox{on } \partial \omega .\\
\end{array}
\right.
$$
By uniqueness of the positive solutions of (\ref{eqn:principaleigenvalueh}), it follows that $\phi = \phi_\Lambda$ (up to multiplication by a positive constant), and
$$\mu_{h,f} (\Lambda)=\Lambda^2 - \Lambda c .$$
Since $\mu_{h,f} (0) <0$ and $\Lambda \geq 0$, it follows that $\Lambda >0$, whence $c\geq c^*$ by definition of $c^*$ (see (\ref{eqn:c*})).

\subsection{The left limits for concentration} \label{sec:leftlimits}

Lastly, we show the convergence of $Y$ as $x\rightarrow -\infty$ to a constant $Y_\infty \in (0,1)$. We have already shown in Section \ref{sec:Tbounded} the existence of such a constant in $[0,1)$. Let us prove that $Y_\infty >0$. We argue by contradiction and assume that $Y_\infty =0$. First, since $c \geq c^*$, there exists $\lambda >0$ such that
$$
\mu_{h,f} (\lambda ) =\lambda ^2  -c\lambda.
$$
Besides, we then have that
\begin{equation}\label{eqn:lambda1}
 \lambda ^2  -c\lambda = \mu_{h,f} (\lambda ) \leq \mu_{h,0} (\lambda ) - \min_{y \in \overline{\omega}} \frac{\partial f}{\partial T} (y,0) .
\end{equation}
Since $T$ is bounded, there exists a constant $C_0 >0$ such that $T(x,y) \leq C_0 e^{-\lambda x}$ for all $x \leq 0$ and $y \in \overline{\omega}$. We then show that there exists $\gamma , \delta \geq 0$ such that $T(x,y) \leq \gamma e^{\delta x}$ for $x \leq 0$.

Indeed, let $$\varepsilon = \min \left( \frac{\mu_{h,0} (0)}{2} , \frac{1}{2} \min_{y \in \overline{\omega}} \frac{\partial f}{\partial T} (y,0) \right) >0,$$ and $A\geq 0$ such that
$$\forall x \leq -A, \ \forall y \in \overline{\omega} , \ \frac{\partial f}{\partial T} (y,0)Y(x,y) \leq \varepsilon .$$
Such a $A$ exists since $Y(-\infty ,.)=Y_\infty =0$. As a consequence of the continuity of $\mu_{h,f}$ and (\ref{eqn:lambda1}), there exists $\Lambda > \lambda$ such that
\begin{equation}\label{eqn:Lambda1}
-\mu_{h,0} (\Lambda ) - c\Lambda + \Lambda ^2 < \frac{1}{2} \min_{y \in \overline{\omega}} \frac{\partial f}{\partial T} (y,0) < -\varepsilon .
\end{equation}
We denote by $U$ the positive function defined by
$$T(x,y)=U(x,y) e^{-\Lambda x} \phi_{0,\Lambda} (y),$$
where $\phi_{0,\Lambda}$ solves (\ref{eqn:principaleigenvalueh}) with the parameter $\Lambda $ and $f=0$, normalized so that $\| \phi_{0,\Lambda} \|_{L^2 (\omega )} =1$. One has $U(-\infty ,.)=0$ as $T(x,y) \leq C_0 e^{-\lambda x}$ for all $x \leq 0$ and $\lambda < \Lambda$, and $\partial_n U =0$ on $\partial \Omega$. Furthermore, it is straightforward to check that
$$ \Delta U + (c-u(y) - 2\Lambda )U_x +2\frac{\nabla_y \phi_{0,\Lambda} }{\phi_{0,\Lambda } }.\nabla_y U $$
$$+(\Lambda^2 -\mu_{h,0} (\Lambda ) +\frac{\partial h}{\partial T} (y,0) - \frac{h(y,T(x,y))}{T(x,y)} -c\Lambda + g(x,y))U  = 0 \mbox{ in } \Omega ,$$
where
$$g(x,y) = \frac{f(y,T(x,y))}{T(x,y)} Y(x,y) \leq \frac{\partial f}{\partial T} (y,0) Y(x,y) \leq \varepsilon$$
for all $(x,y) \in (-\infty , -A]\times \omega$. Besides, we have that
$$\frac{\partial h}{\partial T} (0,y) - \frac{h(y,T(x,y))}{T(x,y)} \leq 0 \mbox{ in } \Omega . $$ Therefore, we have
$$\Delta U +(c-u(y)-2\Lambda )U_x +2\frac{\nabla_y \phi_{0,\Lambda} }{\phi_{0,\Lambda} }.\nabla_y U + (\Lambda^2 -\mu_{h,0} (\Lambda )  -c\Lambda + \varepsilon )U \geq 0$$
for all $(x,y) \in (-\infty , -A]\times \omega$.

Because of (\ref{eqn:Lambda1}), we shall now apply the maximum principle to the previous operator, and look for a suitable super-solution. Since $\varepsilon \leq \mu_{h,0} (0) /2$, there exists $\delta >0$ such that
$$\delta^2 + c\delta - \mu_{h,0} (-\delta ) + \varepsilon < 0 .$$
One can then check that the function
$$\overline{U} (x,y) = e^{(\Lambda + \delta )x} \times \frac{\phi_{0,-\delta } (y)}{\phi_{0, \Lambda} (y)} ,$$ where $\phi_{0,- \delta}$ solves (\ref{eqn:principaleigenvalueh}) with the parameter $\Lambda $ and $f=0$,
satisfies
$$\Delta \overline{U} + (c-u(y)-2\Lambda ) \overline{U}_x +2 \frac{\nabla_y \phi_{0, \Lambda}}{\phi_{0, \Lambda}} . \nabla_y \overline{U} + (\Lambda ^2 - \mu_{h,0} (\Lambda ) -c \Lambda + \varepsilon )\overline{U}$$
$$=(\delta^2 +c\delta -\mu_{h,0} (-\delta )+\varepsilon )\overline{U} \leq 0 \mbox{ in } \Omega ,$$
and
$$\frac{\partial \overline{U}}{\partial n} =0 \mbox{ on } \partial \Omega .$$
It follows from the maximum principle that the difference $\overline{U} -U$ can not attain an interior negative minimum. Moreover, $\overline{U} \geq 0$ and one can normalize the function $\phi_{0,-\delta}$ so that $U(-A,y) \leq \overline{U} (-A ,y)$ for all $y \in \overline{\omega}$. Finally, both $U$ and $\overline{U}$ tend to $0$ as $x\rightarrow -\infty$. We conclude that
$$\forall x \leq -A, \ \forall y \in \overline{\omega}, \ U(x,y)\leq \overline{U} (x,y) . $$
In other words,
$$\forall x \leq -A, \ \forall y \in \overline{\omega}, \ T(x,y)\leq e^{\delta x} \phi_{0,-\delta} (y) \leq \gamma_1 e^{\delta x},$$
where $\gamma_1 = \max_{y\in \overline{\omega}} \phi_{0, -\delta} (y)$. Since $T$ is bounded, we also have
$$\forall x \in (-A,0], \forall y \in \overline{\omega}, \ T(x,y) \leq \gamma_2 e^{\delta x} ,$$
where $\gamma_2 = \|T \|_{L^\infty (\Omega)} e^{\delta A}$. Whence with $\gamma = \max (\gamma_1 , \gamma_2 )$ and for $x \leq 0$ we have
\begin{equation}\label{eqn:majT}
T(x,y) \leq \gamma e^{\delta x} .
\end{equation}

We now claim that
\begin{equation}\label{eqn:Lambda2}
M := \limsup_{x\rightarrow -\infty , \ y\in \overline{\omega}} \frac{Y_x (x,y)}{Y(x,y)} =0 .
\end{equation}
From the Harnack inequality and the fact that $f(y,T)$ is bounded, we know that $|\nabla Y | /Y$ is globally bounded. Therefore, $M$ is finite. Furthermore, $M \geq 0$ because $Y>0 =Y(-\infty ,.)$. Let now $(x_k ,y_k)_{n\in \mathbb{N}}$ a sequence of points in $\mathbb{R} \times \overline{\omega}$ such that $x_k \rightarrow -\infty$ and
$$\frac{Y_x (x_k ,y_k)}{Y(x_k ,y_k)} \rightarrow M \mbox{ as } k\rightarrow +\infty . $$
Up to extraction of some subsequence, one can assume that $y_k \rightarrow y_\infty \in \overline{\omega}$ as $k \rightarrow +\infty$. Consider now the functions
$$Y_k (x,y) =\frac{Y(x+x_k ,y)}{Y(x_k ,y_k)} .$$
They are locally bounded in $\overline{\Omega}$ and satisfy
$$
\left\{
\begin{array}{rcll}
\displaystyle \mbox{Le}^{-1} \Delta Y_k +(c-u(y))Y_{k,x} & = & f(y,T(x+x_k ,y))Y_k & \mbox{in } \Omega,\\
\displaystyle \frac{\partial Y_k}{\partial n}  & = & 0 & \mbox{on } \partial \Omega.
\end{array}
\right.
$$
Moreover, $f(y,T(x+x_k ,y)) \rightarrow 0$ locally uniformly in $\overline{\Omega}$ as $k\rightarrow +\infty $ because $T(-\infty ,.)=0$ and $f(y,0)=0$ for all $y \in \overline{\omega}$. From standard elliptic estimates, up to extraction of some subsequence, the functions $Y_k$ converge weakly in $W_{loc}^{2,p} (\overline{\Omega})$ for all $1\leq p< +\infty$ and strongly in $C_{loc}^{1,\beta } (\overline{\Omega})$ for $0 \leq \beta <1$ to a solution $Z$ of
$$
\left\{
\begin{array}{rcll}
\displaystyle \mbox{Le}^{-1} \Delta Z + (c-u(y)) Z_x & = & 0 & \mbox{in } \Omega ,\\
\displaystyle \frac{\partial Z}{\partial n}  & = & 0 & \mbox{on } \partial \Omega .
\end{array}
\right.
$$
Furthermore, $Z(0, y_\infty ) =1$, $Z\geq 0$ and thus $Z>0$ in $\overline{\Omega}$ from the strong maximum principle and the Hopf lemma. We also have that $Z_x /Z \leq M $ in $\overline{\Omega}$ and $Z_x (0,y_\infty ) /Z(0,y_\infty ) = M$, owing to the definition of $M$ and of the sequence $(x_k ,y_k )$. However, the function $W(x,y)=~Z_x (x,y) / Z(x,y)$ satisfies the equation
$$
\left\{
\begin{array}{rcll}
\displaystyle \mbox{Le}^{-1} \Delta W + 2\mbox{Le}^{-1} \frac{\nabla Z}{Z} . \nabla W +(c-u(y)) W_x & = & 0 & \mbox{in } \Omega ,\\
\displaystyle \frac{\partial W}{\partial n}  & = & 0 & \mbox{on } \partial \Omega .
\end{array}
\right.
$$
Therefore, by the maximum principle and the Hopf lemma, $W(x,y)=M$ for all $(x,y) \in \overline{\Omega}$. In other words, $Z(x,y) =e^{M x} \psi (y)$ where $\psi$ positive function in $\overline{\omega}$ and verifies
$$
\left\{
\begin{array}{rcll}
\displaystyle \mbox{Le}^{-1} \Delta \psi + \mbox{Le}^{-1} M^2 \psi + M (c-u(y)) \psi & = & 0 & \mbox{in }\omega ,\\
\displaystyle \frac{\partial \psi}{\partial n}  & = & 0 & \mbox{on } \partial \omega .
\end{array}
\right.
$$
As a consequence, by uniqueness of the principal eigenvalue  of (\ref{eqn:principaleigenvalue}),
$$M^2 +c M \mbox{Le}= \nu (-M \mbox{Le}).$$
The left-hand side is nonnegative (recall that $c$ is positive, as shown in Section \ref{sec:Tbounded}) while the right-hand side is nonpositive (recall that $\nu$ concave with $\nu (0) = \nu ' (0) =0$). As a conclusion, $M =0$ and then $Z=\psi$ principal eigenfunction of (\ref{eqn:principaleigenvalue}) with parameter $0$ and $\psi (y_\infty) =1$, namely $\psi =1$ in $\overline{\omega}$. Thus, $Z=1$ in $\overline{\Omega}$.
\\

Fix now $\beta >0$ such that $\beta < \delta$ with $\beta$ as in (\ref{eqn:majT}). It follows from (\ref{eqn:Lambda2}) that there exists $A \geq 0$ such that $Y_x (x,y)/Y(x,y) \leq \beta $ for all $x \leq -A$ and $y \in \overline{\omega}$. It follows immediately that
\begin{equation}\label{eqn:minY}
\forall x\leq -A,\ \forall y \in \overline{\omega} , \ Y(x,y) \geq \kappa e^{\beta x} ,
\end{equation}
where $\kappa = e^{-\beta A} \times \min_{y \in \overline{\omega}} Y(-A ,y)>0$.

As we have shown in the proof of (\ref{eqn:Lambda2}), there exists a sequence $(x_k ,y_k)_{k \in \mathbb{N}}$ such that $x_k \rightarrow -\infty$ and the functions $(x,y)\rightarrow Y(x+x_k ,y)/Y(x_k ,y_k)$ converge to the constant 1 at least in $C_{loc}^1 (\overline{\Omega})$ as $k\rightarrow +\infty$. Without loss of generality, one can assume that $x_k \leq -A \leq 0$ for all $k$. Now use the fact that $Y_\infty = Y(-\infty ,.)=Y_x (-\infty ,.)=0$ and integrate equation (\ref{eqn:sysfront}) satisfied by $Y$ over the domain $(-\infty , x_k) \times \omega$. One obtain
\begin{equation}\label{eqn:eqint2}
\mbox{Le}^{-1} \int_{\omega} Y_x (x_k ,y)dy + c\int_\omega Y(x_k ,y)dy - \int_\omega u(y)Y(x_k ,y)dy \leq \gamma |\omega | \delta^{-1} e^{\delta x_k} \times \max_{y\in \overline{\omega}} \frac{\partial f}{\partial T} (y,0),
\end{equation}
because of (\ref{eqn:majT}) and since $f(T)Y \leq \frac{\partial f}{\partial T} (y,0) T$. Furthermore, as $Y(x+x_k ,y)/Y(x_k ,y_k) \rightarrow 1$ in $C_{loc}^1 (\overline{\Omega})$ and since $u(y)$ is bounded in $\omega$ and has mean zero, it follows that
$$\displaystyle \frac{\int_\omega Y_x (x_k ,y)dy}{\int_\omega Y(x_k ,y) dy} \rightarrow 0$$
and
$$\displaystyle \frac{\int_\omega u(y) Y (x_k ,y)dy}{\int_\omega Y(x_k ,y) dy} \rightarrow | \omega |^{-1} \int_\omega u(y) dy =0$$
as $k\rightarrow +\infty$. Putting that together with (\ref{eqn:eqint2}), one gets that
$$\frac{c}{2} \int_\omega Y(x_k ,y)dy \leq \gamma | \omega | \delta^{-1} e^{\delta x_k} \times \max_{y\in \overline{\omega}} \frac{\partial f}{\partial T} (y,0) $$
for $k$ large enough, because $c>0$. But (\ref{eqn:minY}) together with $x_k \leq -A$ then yields
$$\frac{c \kappa | \omega |}{2}e^{\beta x_k} \leq  \gamma | \omega | \delta^{-1} e^{\delta x_k}\times  \max_{y\in \overline{\omega}} \frac{\partial f}{\partial T} (y,0)$$
for $k$ large enough. Since $0 < \beta < \delta$, one gets a contradiction by passing to the limit $x_k \rightarrow -\infty$.

As a conclusion, $Y_\infty =0$ is impossible. Therefore, $Y_\infty \in (0,1)$ and the proof of Theorem~\ref{th:thqualit} is achieved.~$\Box$

\section{Existence of fronts with non-minimal speeds} \label{sec:exist}

Here, we prove Part (a) of Theorem \ref{th:thexist}.  We assume that $\mu_{h,f} (0) < 0$ and we let $c > \max (0,c^* )$. First, we will construct sub and super-solutions of (\ref{eqn:sysfront}). We will then use a fix point theorem on bounded cylinders to construct approximate solutions. Lastly, by passing to the limit of an infinite cylinder, we will obtain a solution of (\ref{eqn:sysfront}) with the wanted qualitative properties. This now standard procedure has already been applied to show the existence of fronts in \cite{hamel-quenching,berestycki3,hamel-adiabatic} .

\subsection{Sub- and supersolutions in $\overline{\Omega}$}

Note first that the constant 1 is a super-solution for $Y$.

\subsubsection*{Supersolution for T}
We then construct a supersolution for the $T$-equation (\ref{eqn:sysfront}) with $Y=1$.
Since $c> c^*$, let $\lambda_c$ the smallest nonnegative root of $k(\lambda ) =c \lambda$. The real number $\lambda _c$ is in fact positive thanks to $k(0) = - \mu_{h,f} (0) >0$. Then let $\overline{T}$ be the function defined in $\overline{\Omega}$ by
$$\overline{T} (x,y)=\phi_{\lambda_c} (y) e^{-\lambda_c x} >0 .$$
Here $\phi_{\lambda_c}$ is the positive principal eigenfunction of (\ref{eqn:principaleigenvalueh}) with $\lambda =\lambda_c$, normalized so that $\| \phi_{\lambda_c} \|_{L^{\infty} (\omega )} =1$. The function $\overline{T}$ satisfies the Neumann boundary conditions on $\partial \Omega$, and is a super-solution for the equation on $T$ in (\ref{eqn:sysfront}) with $Y=1$, i.e
$$\Delta \overline{T} + (c-u(y)) \overline{T}_x +f(y,\overline{T}) -h(y,\overline{T})$$
$$\leq \Delta \overline{T} + (c-u(y)) \overline{T}_x +(\frac{\partial f}{\partial T} (y,0)-\frac{\partial h}{\partial T} (y,0)) \overline{T} =0 \mbox{ in } \overline{\Omega}.$$

\subsubsection*{Sub-solution for Y}\label{sec:subsolY}

Since $\nu (0) = \nu '(0) =0 < c$, one can choose $\beta >0$ small enough so that
\begin{equation}\label{eqn:beta}
\left\{
\begin{array}{l}
0<\beta < \lambda_c ,\\
\nu (\beta \mbox{Le})-\beta^2 +c\beta \mbox{Le} >0 ,
\end{array}
\right.
\end{equation}
and $\gamma > 0$ large enough so that
\begin{equation}\label{eqn:gamma}
\left\{
\begin{array}{l}
\displaystyle \gamma \times \min_{\overline{\omega}} \psi_{\beta \text{Le}} \geq 1 ,\\
\displaystyle \gamma \mbox{Le}^{-1} (\nu (\beta \mbox{Le})-\beta^2 +c\beta \mbox{Le}) \times \min_{\overline{\omega}} \psi_{\beta \text{Le}} > \max_{y \in \overline{\omega}} \frac{\partial f}{\partial T} (y,0),
\end{array}
\right.
\end{equation}
where $\psi_{\beta \text{Le}}$ is the positive eigenfunction of (\ref{eqn:principaleigenvalue}) with $\lambda=\beta \mbox{Le}$, normalized in such a way that $\|\psi_{\beta \text{Le}}\|_{L^\infty (\omega )}~=~1$. Let $\underline{Y}$ be defined by
$$\underline{Y} (x,y) = \max(0,1-\gamma \psi_{\beta \text{Le}} (y) e^{-\beta x}).$$
Note that $\underline{Y} =0$ for $x \leq 0$. Let us check that $\underline{Y}$ is a sub-solution for (\ref{eqn:sysfront}) with $T=\overline{T}$. Note first that $\underline{Y}$ satisfies the Neumann boundary conditions on $\partial \Omega$. Moreover, when $\underline{Y} >0$, then $x>0$ and
\begin{eqnarray*}
&& \text{Le}^{-1} \Delta \underline{Y} + (c-u(y))\underline{Y}_x - f(y,\overline{T})\underline{Y}\\
&&  \geq \gamma \text{Le}^{-1} ( \nu(\beta \text{Le}) -\beta^2 +c\beta \text{Le}) \psi_{\beta \text{Le}} (y) e^{-\beta x} - \frac{\partial f}{\partial T} (y,0) \phi_{\lambda_c} (y) e^{-\lambda_c x} (1 - \gamma \psi_{\beta \text{Le}} (y) e^{-\beta x})\\
&& \geq \gamma \text{Le}^{-1} ( \nu(\beta \text{Le}) -\beta^2 +c\beta \text{Le}) \psi_{\beta \text{Le}} (y) e^{-\beta x} - \frac{\partial f}{\partial T} (y,0) e^{-\beta x} \geq 0 ,
\end{eqnarray*}since $f$ of the KPP-type, $0<\phi_{\lambda_c} (y) \leq 1$ in $\overline{\omega}$ and because of (\ref{eqn:beta})-(\ref{eqn:gamma}).

\subsubsection*{Sub-solution for T}\label{sec:subsolT}

Lastly, we will construct a sub-solution for $T$ with $Y=\underline{Y}$. Recall that $k(\lambda_c )=c\lambda_c $. We first show that $k'(\lambda_c ) < c$. Indeed, since $k(0)>0$ and $\lambda_c$ is the smallest positive root of $k(\lambda )=c\lambda$, we have $k'(\lambda_c )\leq c$. Furthermore, if $k'(\lambda_c) =c$, then $k(\lambda )\geq c\lambda $ for all $\lambda \in \mathbb{R}$ by convexity of $k$, whence $c^* \geq c$, which is impossible. We conclude, as announced, that $k' (\lambda_c ) < c$.

The above allows us to choose $\eta >0$ small enough so that
\begin{equation}\label{eqn:eta}
\left\{
\begin{array}{l}
0<\eta < \min (\beta ,\alpha \lambda_c) ,\\
\varepsilon :=c(\lambda_c +\eta )-k(\lambda_c + \eta )>0 ,
\end{array}
\right.
\end{equation}
where $\alpha >0$ such that $f(y,.)$ and $h(y,.)$ are of class $C^{1,\alpha } ([0,s_0 ])$ for some $s_0 >0$ uniformly in $y \in \overline{\omega}$. Let $M \geq 0$ such that
\begin{equation}\label{eqn:eqM}
\left\{
\begin{array}{l}
\displaystyle f(y,s)\geq \frac{\partial f}{\partial T} (y,0)s - Ms^{1+\alpha } ,\\
\displaystyle h(y,s)\leq \frac{\partial h}{\partial T}(y,0)s + Ms^{1+\alpha } ,\\
\end{array}
\right. \mbox{ for all } s \in [0,s_0]  \mbox{ and for all } y \in \overline{\omega}.
\end{equation}
Now take $x_0 \geq 0$ sufficiently large so that
$$\underline{Y} (x,y) = 1 -\gamma \psi_{\beta \text{Le}} (y) e^{-\beta x} \mbox{ for all } (x,y) \in (x_0 ,+\infty ) \times \overline{\omega} .$$
Next, let $\delta >0$ large enough so that
\begin{equation}\label{eqn:delta1}
\left\{
\begin{array}{l}
\displaystyle \phi_{\lambda_c} (y) e^{-\lambda_c x} - \delta \phi_{\lambda_c + \eta} (y) e^{-(\lambda_c +\eta ) x} \leq s_0 \mbox{ in } \overline{\Omega} ,\\
\displaystyle \phi_{\lambda_c} (y) e^{-\lambda_c x} - \delta \phi_{\lambda_c + \eta} (y) e^{-(\lambda_c +\eta ) x} \leq 0 \mbox{ \ in } (-\infty ,x_0 ] \times \overline{\omega} ,\\
\displaystyle \delta \varepsilon \times \min_{\overline{\omega}} \phi_{\lambda_c + \eta } \geq \gamma \max_{y \in \overline{\omega}} \frac{\partial f}{\partial T} (y,0)+ 2M .
\end{array}
\right.
\end{equation}
Finally, we define, for all $(x,y) \in \overline{\Omega}$,
$$\underline{T} (x,y) = \max \left( 0 , \phi_{\lambda _c} (y) e^{-\lambda_c x} -\delta \phi_{\lambda_c + \eta} (y) e^{-(\lambda_c +\eta )x} \right) .$$
The function $\underline{T}$ satisfies the Neumann boundary conditions on $\partial \Omega$. Let us now check that $\underline{T}$ is a sub-solution to (\ref{eqn:sysfront}) with $Y=\underline{Y}$. Note first that $0 \leq \underline{T} \leq s_0$ in $\overline{\Omega}$. Moreover, if $\underline{T} (x,y) > 0$, then $x>x_0 \geq 0$ whence $0\leq \underline{Y}(x,y)=1-\gamma \psi_{\beta \text{Le}} (y) e^{-\beta x}$.
Then, in that case, we have:
\begin{eqnarray*}
&& \Delta \underline{T} + (c-u(y))\underline{T}_x + f(y,\underline{T})\underline{Y} - h(y,\underline{T})\\
&& \geq \Delta \underline{T} + (c-u(y))\underline{T}_x + (\frac{\partial f}{\partial T} (y,0)\underline{T}-M\underline{T}^{1+\alpha})(1-\gamma \psi_{\beta \text{Le}} (y) e^{-\beta x})- \frac{\partial h}{\partial T} (y, 0)\underline{T}-M\underline{T}^{1+\alpha}\\
&& \geq -\delta ( k(\lambda_c +\eta ) -c(\lambda_c +\eta))\phi_{\lambda_c +\eta} (y) e^{-(\lambda_c +\eta)x} -\frac{\partial f}{\partial T} (y,0)\gamma \underline{T}\psi_{\beta \text{Le}} (y) e^{-\beta x} -2M\underline{T}^{1+\alpha}\\
&& \geq \delta \varepsilon \phi_{\lambda_c +\eta} (y) e^{-(\lambda_c +\eta)x} - \frac{\partial f}{\partial T} (y,0)\gamma e^{-(\lambda_c + \beta )x} - 2Me^{-\lambda_c (1+\alpha )x}\\
&& \geq  (\delta \varepsilon \phi_{\lambda_c +\eta} (y) -\frac{\partial f}{\partial T} (y,0)\gamma - 2M)e^{-(\lambda_c + \eta )x} \geq 0 ,
\end{eqnarray*}
because of (\ref{eqn:eta}), (\ref{eqn:eqM}), (\ref{eqn:delta1}) and since $0<\phi_{\lambda_c +\eta} (y)$, $0 < \psi_{\beta \text{Le}} (y) \leq 1$ in $\overline{\omega}$.
\subsection{The finite cylinder problem}

Here, we construct a solution of (\ref{eqn:sysfront}) in a finite cylinder $\Omega_a =(-a,a)\times \omega$ with $a>0$. Let $C(\overline{\Omega_a})$ denote the space of continuous functions in $\overline{\Omega_a}$, with the usual sup-norm. Observe that $0 \leq \underline{T} \leq \overline{T}$ and $0\leq \underline{Y} \leq 1$ in $\overline{\Omega}$. We denote by $E_a$ the set
$$E_a =\{(T,Y) \in C(\overline{\Omega_a};\mathbb{R}^2), \ \underline{T}\leq T \leq \overline{T} \mbox{ and } \underline{Y} \leq Y \leq 1 \mbox{ in } \overline{\Omega_a}\}.$$
The set $E_a$ is a convex closed bounded subset of the Banach space $C(\overline{\Omega_a};\mathbb{R}^2)$.

We now consider a fixed point problem for an approximation of the travelling front solution in $\Omega_a$. For any pair $(T_0 ,Y_0) \in E_a$, let $(T,Y)=\Phi_a (T_0 ,Y_0 )$ be the unique solution of
\begin{equation*}
\left\{
\begin{array}{rcll}
\Delta T +(c-u(y))T_x - K_a T& = & -f(y,T_0)Y_0  + h(y,T_0) - K_a T_0& \mbox{in } \Omega_a ,\\
\text{Le}^{-1} \Delta Y + (c-u(y)) Y_x -f(y,T_0)Y & = & 0 & \mbox{in } \Omega_a ,\\
\end{array}
\right.
\end{equation*}
together with the boundary conditions
\begin{equation*}
\left\{
\begin{array}{rl}
\displaystyle T(\pm a,y) = \underline{T}(\pm a,y), \ Y(\pm a,y)= \underline{Y} (\pm a,y) & \mbox{for } y \in \overline{\omega} ,\\
\displaystyle \frac{\partial T}{\partial n} = \frac{\partial Y}{\partial n} = 0 & \mbox{on } [-a ,a] \times \partial \omega .\\
\end{array}
\right.
\end{equation*}
Since $h$ is in $C^1 (\overline{\omega} \times [0,+\infty);\mathbb{R})$, we can assume that $K_a$ is positive and such that for all $y \in \overline{\omega}$
\begin{equation}\label{eqn:Kdec}
s \in [0, \sup_{\Omega_a} \overline{T} \ ] \rightarrow h(y,s) - K_a s \mbox{  \ \   } \mbox{is decreasing}.
\end{equation}
Such a solution $(T,Y)$ exists, belongs to $C(\overline{\Omega_a} ; \mathbb{R}^2)$ and it is unique (see \cite{berestycki2,berestycki4}). To show that the map $\Phi_a$ has a fixed point, we will show that the set $E_a$ is invariant by $\Phi_a$, and that the map $\Phi_a$ is compact.

\subsubsection*{$E_a$ is invariant by $\Phi_a$}

Let $(T_0 , Y_0)$ be any element of $E_a$, and $(T,Y)=\Phi_a (T_0 ,Y_0)$. One can check that $\underline{T}$ satisfies
$$\Delta \underline{T} + (c-u(y))\underline{T}_x - K_a \underline{T} \geq -f(y,\underline{T})\underline{Y} + h(y,\underline{T}) - K_a \underline{T} \geq -f(y,T_0)Y_0 + h(y, T_0) - K_a T_0 , $$
where the last inequality follows from (\ref{eqn:Kdec}) and the monotonicity of $f$. Furthermore, $\underline{T}$ satisfies the same boundary conditions as $T$ on the boundary of $\Omega_a$. The weak maximum principle implies that $\underline{T} \leq T$ in $\overline{\Omega_a}$. The inequalities $T \leq \overline{T}$, $\underline{Y}
 \leq Y$ and $Y\leq 1$ in $\overline{\Omega_a}$ can be checked similarly.

We conclude that $\Phi_a$ leaves $E_a$ invariant.

\subsubsection*{The map $\Phi_a$ is compact}

We introduce $(k_1 ,j_1) =\Phi_a (\overline{T},1)$ and $(k_2 ,j_2)=\Phi_a (\underline{T},1)$. For any pair $(T_0 ,Y_0) \in E_a$ and $(T,Y)=\Phi_a (T_0 ,Y_0)$, one has
$$\Delta k_1 + (c-u(y))k_{1,x} - K_a k_1 = -f(y,\overline{T}) + h(y,\overline{T}) -K_a \overline{T} \leq -f(y,T_0) Y_0 + h(y,T_0) - K_a T_0 \mbox{ in } \Omega_a ,$$
and thus $T \leq k_1$ in $\overline{\Omega_a}$. Similarly, we have
$$\text{Le}^{-1} \Delta j_2 + (c-u(y))j_{2,x} - f(y,T_0)j_2 = (f(y,\underline{T})-f(y,T_0))j_2 \leq 0 \mbox{ in } \Omega_a ,$$
and thus $Y \leq j_2$ in $\overline{\Omega_a}$. Therefore, we obtain
\begin{equation}\label{eqn:inequalities1}
\left\{
\begin{array}{l}
\underline{T} \leq T \leq k_1 \leq \overline{T}\\
\underline{Y} \leq Y \leq j_2 \leq 1 \\
\end{array}
\right.
\mbox{ in } \overline{\Omega_a} ,
\end{equation}
for any pair $(T_0 ,Y_0) \in E_a$ and $(T,Y)=\Phi_a (T_0 ,Y_0)$.

Let $(T_0^n ,Y_0^n)$ be a sequence in $E_a$ and
$$(T^n ,Y^n)=\Phi_a (T_0^n ,Y_0^n ).$$
By standard elliptic estimates up to the boundary, the sequence $(T^n ,Y^n)$ is bounded in $C^1 (D; \mathbb{R}^2)$ norm, for any compact subset
$$D \subset \Sigma_a = \overline{\Omega_a} \setminus \{ \pm a \} \times \partial \omega . $$
Therefore, using the diagonal extraction process, there exists a subsequence, still denoted by $(T^n ,Y^n )$ which converges locally uniformly in $\Sigma_a$ to a pair $(T,Y)$ of continuous functions in $\Sigma_a$. Since each $(T^n ,Y^n)$ satisfies (\ref{eqn:inequalities1}) in $\overline{\Omega_a}$, it follows that $(T,Y)$ satisfies (\ref{eqn:inequalities1}) in $\Sigma_a$. Furthermore, as we have
$$
\left\{
\begin{array}{l}
k_1 (\pm a ,y)=\underline{T}(\pm a ,y)\\
j_2 (\pm a ,y)=\underline{Y}(\pm a ,y)\\
\end{array}
\right.
\mbox{ in } \overline{\omega},
$$
and since $\underline{T}$, $\underline{Y}$, $k_1$ and $j_2$ are continuous in $\overline{\Omega_a}$, the functions $(T,Y)$ can be extended in $\overline{\Omega_a}$ by two continuous functions, still denoted by $(T,Y)$, satisfying (\ref{eqn:inequalities1}) in $\overline{\Omega_a}$.

For any $\varepsilon >0$, there exists $\kappa >0$ such that
$$
\left\{
\begin{array}{l}
0 \leq k_1 - \underline{T} \leq \varepsilon\\
0 \leq j_2 - \underline{Y} \leq \varepsilon\\
\end{array}
\right.
\mbox{ in } [-a,-a+\kappa ] \times \overline{\omega} \ \cup \ [a-\kappa ,a]\times \overline{\omega} ,
$$
and thus $| T^n - T| \leq \varepsilon$ and $|Y^n -T| \leq \varepsilon$ in the same sets, for all $n$. On the other hand, the sequence $(T^n ,Y^n )$ converges uniformly to $(T,Y)$ in $[-a+\kappa , a-\kappa ]\times \overline{\omega}$. Therefore, $(T^n ,Y^n)$ converges uniformly to $(T,Y)$ in $\overline{\Omega_a}=[-a ,a]\times \overline{\omega}$ and thus the map $\Phi_a$ is compact.

\subsubsection*{A fixed point of $\Phi_a$}

One then concludes from the Schauder fixed point theorem that $\Phi_a$ has a fixed point in $E_a$. In other words, there exists a solution $(T_a ,Y_a) \in E_a$ of the problem
\begin{equation}\label{eqn:finitecylinder}
\left\{
\begin{array}{rcll}
\Delta T_a + (c-u(y))T_{a,x}+f(y,T_a) Y_a -h(y,T_a) & = & 0 & \mbox{ in } \Omega_a ,\\
\text{Le}^{-1} \Delta Y_a + (c-u(y))Y_{a,x} -f(y,T_a) Y_a  & = & 0 & \mbox{ in } \Omega_a ,\\
\end{array}
\right.
\end{equation}
with the boundary conditions
\begin{equation}\label{eqn:finitecylinderboundary}
\left\{
\begin{array}{rl}
\displaystyle T_a(\pm a,y) = \underline{T}(\pm a,y), \ Y_a(\pm a,y)= \underline{Y} (\pm a,y) & \mbox{for } y \in \overline{\omega} , \\
\displaystyle \frac{\partial T_a}{\partial n} = \frac{\partial Y_a}{\partial n} = 0 & \mbox{on } [-a ,a] \times \partial \omega . \\
\end{array}
\right.
\end{equation}
Furthermore, we have $0 \leq \underline{T} \leq T_a \leq \overline{T}$ and $0\leq \underline{Y} \leq Y_a \leq 1$ in $[-a ,a] \times \overline{\omega}$.

\subsection{Passage to the infinite cylinder}

Let now $(a_n)_{n\in \mathbb{N}}$ be an increasing sequence of positive numbers such that $a_n \rightarrow +\infty$ as $n \rightarrow + \infty$. Let $(T_{a_n},Y_{a_n})_{n\in \mathbb{N}}$ be a sequence of solutions of (\ref{eqn:finitecylinder})-(\ref{eqn:finitecylinderboundary}) with $a=a_n$. From standard elliptic estimates up to the boundary, the sequence $(T_{a_n},Y_{a_n})$ is bounded in $C_{loc}^{2,\alpha} (\overline{\Omega})$ (remember that the flow u is of class $C^{0,\alpha} (\overline{\omega})$ and $f$, $h$ are locally Lipschitz-continuous). Up to extraction of a subsequence, the functions $(T_{a_n},Y_{a_n})$ then converge in $C_{loc}^2 (\overline{\Omega})$ to a pair $(T,Y) \in C^2 (\overline{\Omega})$ of solutions of
$$
\left\{
\begin{array}{rcll}
\Delta T + (c-u(y))T_x + f(y,T)Y - h(y,T) & = & 0 & \mbox{ in } \Omega ,\\
\text{Le}^{-1} \Delta Y + (c-u(y))T_x - f(y,T)Y & = & 0 & \mbox{ in } \Omega ,\\
\end{array}
\right.
$$
with the Neumann boundary conditions
$$\frac{\partial T}{\partial n} = \frac{\partial Y}{\partial n} = 0 \mbox{ on } \partial \Omega ,$$
and
$$0\leq \underline{T} \leq T \leq \overline{T} \mbox{ ,  } 0\leq \underline{Y} \leq Y \leq 1 \mbox{  in } \overline{\Omega} .$$
In particular, we have $T(+\infty ,y)=0$ and $Y(+\infty ,y) =1$ uniformly in $y\in \overline{\omega}$. Furthermore, since $\underline{Y} (x,y)$ and $\underline{T} (x,y)$ are positive for large $x$, the strong maximum principle implies that $Y>0$ and $T>0$ in $\overline{\Omega}$. Moreover, since $f(y,T)>0$, the function $Y$ cannot be identically equal to 1,  whence $Y<1$ in $\overline{\Omega}$ from the strong maximum principle.

It now remains to be shown that $T$ is bounded, and that the functions $(T,Y)$ satisfy the right conditions at $-\infty$.

\subsection{Boundedness of $T$}

Assume for the sake of a contradiction that $T$ is not in $L^\infty (\Omega)$. Since $0 \leq \underline{T} \leq T \leq \overline{T}$ in $\overline{\Omega}$, the only possibility for the function $T$ to grow is on the left. Thus there exists a sequence $(x_n ,y_n)_{n \in \mathbb N}$ in $\mathbb{R} \times \overline{\omega}$ such that
$$T(x_n ,y_n) \rightarrow +\infty \mbox{ and } x_n \rightarrow -\infty$$
as $n\rightarrow +\infty$. Since the function $| \nabla T | / T$ is globally bounded from standard elliptic estimates and the Harnack inequality up to the boundary, it follows that for each $R >0$,
$$\min_{(x,y) \in [x_n -R ,x_n +R] \times \overline{\omega}} T (x ,y) \rightarrow +\infty$$
as $n \rightarrow + \infty$. Let also $m=\min_{y \in \overline{\omega}} f(y,1) >0$. We recall that the function $\nu$ defined in (\ref{eqn:principaleigenvalue}) is concave and that $\nu (0) =0$. Therefore, there exist exactly two real numbers $\rho_\pm$ such that $\rho_{-} < 0 < \rho_{+}$ and
$$\text{Le}^{-1} \nu (-\rho_{\pm} \text{Le}) = \text{Le}^{-1} \rho_{\pm}^2 + c\rho_{\pm} -m .$$
We denote by $\psi_{\pm}$ the two principal eigenfunctions of the problem (\ref{eqn:principaleigenvalue}) with the values $\lambda = - \rho_{\pm} \text{Le}$, normalized so that, say, $\min_{\overline{\omega}} \psi_{\pm} = 1$. The functions $u_{\pm} (x,y) = e^{\rho_{\pm} x} \psi_{\pm} (y)$ then satisfy
$$
\left\{
\begin{array}{rcll}
\displaystyle \text{Le}^{-1} \Delta u_{\pm} + (c-u(y)) u_{\pm ,x} - m u_{\pm} & = & 0 & \mbox{in } \Omega ,\\
\displaystyle \frac{\partial u_{\pm}}{\partial n} & = & 0 & \mbox{on } \partial \Omega .\\
\end{array}
\right.
$$
Fix now any $R >0$ and choose $N \in \mathbb{N}$ so that
$$\min_{(x,y) \in [x_n -R ,x_n +R] \times \overline{\omega}} T (x ,y) \geq 1$$
for all $n \geq N$. Then, as the function $f(y,T)$ is increasing in the variable $T$, we have that $f(y,T) \geq f(y,1) \geq m$ in $[x_n -R ,x_n +R] \times \overline{\omega}$ for all $n \geq N$. Whence, on the same domain,
$$\text{Le}^{-1} \Delta Y + (c-u(y)) Y_x - mY \geq 0 .$$
The function $Y$ also satisfies the Neumann boundary conditions on $\partial \Omega $. Furthermore, $Y \leq 1$ in $\Omega$. It then follows from the weak maximum principle that
$$\forall (x,y) \in [x_n -R ,x_n +R] \times \overline{\omega} \mbox{, } Y(x,y) \leq e^{\rho_+ (x-x_n -R)} \psi_+ (y) + e^{\rho_- (x-x_n +R)} \psi_- (y).$$
Therefore, along the section $x = x_n$, the function $Y$ is small:
$$\limsup_{n\rightarrow +\infty} \max_{y \in \overline{\omega}} Y (x_n ,y) \leq \max \left( \max_{\overline{\omega}} \psi_+ , \max_{\overline{\omega}} \psi_- \right) \times ( e^{- \rho_+ R} + e^{\rho_- R}) .$$
Since $R>0$ can be chosen arbitrary, one concludes that $Y(x_n ,.) \rightarrow 0$ uniformly in $\overline{\omega}$ as $n \rightarrow +\infty$.

Let now $\epsilon > 0$ be any positive real number, and $N \in \mathbb{N}$ such that $Y(x_n ,y) \leq \epsilon$ for all $n \geq N$ and $y \in \overline{\omega}$. Since the function $Y$ satisfies
$$\mbox{Le}^{-1} \Delta Y + (c-u(y))Y_x = f(y,T) Y \geq 0 ,$$
it follows from the weak maximum principle that
$$Y(x,y) \leq \epsilon$$
for all $(x,y) \in [x_n ,x_N] \times \overline{\omega}$ and $n \geq N$ such that $x_n \leq x_N$. Since $x_n \rightarrow -\infty$ as $n \rightarrow +\infty$, one concludes that $Y \leq \epsilon$ in $(-\infty , x_N ] \times \overline{\omega}$. As $Y \geq 0$, we finally obtain that $Y(-\infty,.) =0$ uniformly in $y \in \overline{\omega}$.

\par
We now use the same arguments as in Section \ref{sec:leftlimits}. We have just shown that $Y(-\infty,.) =0$. Furthermore, since $T\leq \overline{T}$, we know that there exist $C_0 >0$ and $\lambda >0$ solution of $k(\lambda )= c\lambda$ such that $T \leq C_0 e^{-\lambda x}$. As already shown in Section \ref{sec:leftlimits}, it then implies that there exist $A$, $\gamma$, $\delta \geq 0$ such that $T \leq \gamma e^{\delta x}$ for all $x \leq -A$, which is in contradiction with $T(x_n,y_n)\rightarrow +\infty$ and $x_n \rightarrow -\infty$ as $n \rightarrow + \infty$.

We conclude that $T$ is bounded.

\subsection{Behavior of the solution on the left}\label{sec:behavleft}

It now only remains to show $T_x (-\infty, .) =Y_x (-\infty,.)=0$. In fact, we show that $T$ and $Y$ converge to constants as $x \rightarrow -\infty$. We will then conclude by standard elliptic estimates.

Since $T$ and $Y$ are globally bounded, standard elliptic estimates and Harnack inequality imply that $\nabla T$ and $\nabla Y$ are globally bounded as well. As in Section \ref{sec:Tbounded}, by integrating the equation (\ref{eqn:sysfront}) satisfied by $Y$ over the domain $(-N,N) \times \omega$ where $N>0$, we obtain
\begin{equation}\label{eqn:eqint3}
\begin{array}{c}
\int_\omega [ \text{Le}^{-1} (Y_x (N,y) - Y_x (-N,y)) + (c-u(y))(Y(N,y)-Y(-N,y))]dy\\
\\
=\int_{(-N,N) \times \omega} f(y,T(x,y))Y(x,y) dx dy .
\end{array}
\end{equation}
The left-hand side is then bounded independently of $N$, whence $$ \int_{\Omega} f(y,T)Y$$ converges. Next, by multiplying the equation (\ref{eqn:sysfront}) satisfied by $Y$ by $Y$ itself, and integrating over the domain $(-N,N) \times \omega$ for $N>0$, we obtain
$$\int_\omega [ \text{Le}^{-1} (Y_x (N,y) Y(N,y) - Y_x (-N,y))Y(-N,y) + \frac{1}{2} (c-u(y))(Y^2 (N,y) -Y^2 (-N,y))]dy$$
$$=\int_{(-N,N) \times \omega} [f(y,T(x,y))Y^2 (x,y) + \text{Le}^{-1} |\nabla Y | ^2 ]dx dy.$$
The left-hand side is again bounded independently of $N$, and so is the integral
$$0 \leq \int_{(-N,N) \times \omega} f(y,T(x,y))Y^2 (x,y) dx dy .$$
Therefore, we have
\begin{equation}\label{eqn:inegintY}
\int_{\Omega} |\nabla Y(x,y)|^2 dxdy < +\infty .
\end{equation}
Choose now any sequence $(x_k)_{k \in \mathbb{N}} \rightarrow -\infty$ and define the translates
$$Y_k (x,y) = Y(x_k +x,y) .$$
The functions $Y_k$ are bounded in $W_{loc}^{2,p} (\overline{\Omega})$ for all $1 < p< \infty$. Therefore, up to the extraction of a subsequence, the functions $Y_k$ converge in $C_{loc}^1 (\overline{\Omega})$ to a function $Y_\infty $. It follows from (\ref{eqn:inegintY}) that $Y_\infty$ is a constant. We now show that this constant does not depend on the choice of the subsequence. Recall that $Y_x (+\infty ,.)=0$ because $Y$ converges to a constant as $x\rightarrow +\infty$ and from standard elliptic estimates. We set $N=-x_k$ in (\ref{eqn:eqint3}) and pass to the limit $k\rightarrow +\infty$. This leads to
$$\int_\omega (c-u(y)) (1-Y_\infty) dy = \int_{\Omega} f(y,T(x,y))Y(x,y)dxdy ,$$
and, since $u$ has zero average over $\omega$,
$$c | \omega | (1 -Y_\infty ) = \int_{\Omega} f(y,T(x,y))Y(x,y) dxdy .$$
Therefore, $Y_\infty$ does not depend on the sequence $(x_k)_{k\in \mathbb{N}}$. Thus, the limit $Y(-\infty ,.) = Y_\infty$ exists and $Y_x (-\infty ,.) = 0$.

Let us now prove that $T(-\infty ,. )=T_x (-\infty ,.) =0$. We integrate the equation (\ref{eqn:sysfront}) satisfied by $T$ over $(-N,N) \times \omega$ for $N>0$. We obtain
$$\int_\omega [(T_x (N,y) - T_x (-N,y)) + (c-u(y))(T(N,y) - T(-N,y))]dy $$ $$+ \int_{(-N,N) \times \omega} f(y,T(x,y))Y(x,y) dxdy
=\int_{(-N,N) \times \omega} h(y,T(x,y))dxdy .$$
The left-hand side is bounded independently of $N$ (recall that $T$ and $T_x$ are bounded) and the function $h$ is nonnegative. Therefore, the integral
$$\int_{\Omega} h(y,T(x,y))dxdy$$
converges.

As in Section \ref{sec:leftlimitT}, we then show that
$$\int_{\Omega} \frac{\partial h }{\partial T} (y,0) T(x,y) dxdy \ < +\infty ,$$
$$\int_{\Omega} | \nabla T |^2 dxdy \ < \ +\infty .$$
We then conclude as in Section \ref{sec:leftlimitT} that $T(-\infty,.) =0$ and thus $T_x (-\infty ,.)=0$ by standard elliptic estimates.

The proof of Part (a) of Theorem \ref{th:thexist} is now complete. $\Box$

\section{Existence of fronts with minimal speed} \label{sec:existc*}

This section is dedicated to the proof of Part (b) of Theorem \ref{th:thexist}. Here, we will assume that
\begin{equation}\label{eqn:condlambda}
\sup_{\lambda \in \mathbb{R}}( \mu_{h,f} (\lambda ) - \lambda^2 )< 0 .
\end{equation}
It then follows immediately from the definition of $c^*$ that $c^* >0$.

Before we begin the proof of Part (b) of Theorem \ref{th:thexist}, observe first that if the average of $\frac{\partial h}{\partial T} (y,0) - \frac{\partial f}{\partial T} (y,0)$ over $\omega$ is less than 0, then condition (\ref{eqn:condlambda}) is satisfied. Indeed, for any $\lambda \in \mathbb{R}$, by dividing (\ref{eqn:principaleigenvalueh}) by $\phi_\lambda$ and integrating over $\omega$, it follows that
$$\mu_{h,f} (\lambda) \leq |\omega|^{-1} \int_\omega (\frac{\partial h}{\partial T} (y,0) - \frac{\partial f}{\partial T} (y,0)) dy$$ because of (\ref{eqn:u}).

Let us now compare the condition (\ref{eqn:condlambda}) with the condition $\mu_{h,f} (0) < 0$. As we said in Remark~\ref{rem:minspeed1}, those hypotheses are equivalent in the case $h$ independent of $y$. Otherwise, it depends on the flow $u$. Indeed, let first $h$ be in the form $h(T)=aT$ with $a \in \mathbb{R}^+$ such that $\mu_{h,f} (0) =0$. Such a $h$ exists because, as one can easily check, $\mu_{0,f} (0) < 0$ and $\mu_{aT ,f} (0) = \mu_{0,f} (0) +a$ for all $a\in \mathbb{R}^+$. Furthermore, from Section \ref{sec:intro}, we know that
$$
\mu_{h,f} ' (0) = - \int_\omega u(y) \phi_0^2 (y) dy,
$$
where $\phi_0$ is a solution of
\begin{equation}
\left\{
\begin{array}{rcll}
\displaystyle -\Delta_y \phi_0 + (\frac{\partial h}{\partial T} (y,0)- \frac{\partial f}{\partial T} (y,0)) \phi_0 & = & \mu_{h,f} (0 ) \phi_0 & \mbox{ in } \omega ,\\
\displaystyle \frac{\partial \phi_0 }{\partial n} & = & 0 & \mbox{ on } \partial \omega ,\\
\end{array}
\right.
\end{equation}
with $L^2 (\omega)$ norm equal to 1. Note that $\phi_0$ is independent of $u$. Thus, if $\phi_0$ is not constant, which is equivalent to say that $\frac{\partial h}{\partial T} (y,0)$ is not constant, a suitable choice of $u$ allows us to obtain any value for $\mu_{h,f} ' (0)$. For instance, we can choose $u$ so that $\mu_{h,f} ' (0) >0$, and then there exists $\lambda >0$ such that
$$\mu_{h,f} (\lambda) - \lambda^2 >0 .$$
Besides, let the sequence $(h_n)_{n\in \mathbb{N}}$ defined by $h_n (T) =h(T)- \frac{1}{n} T = (a-\frac{1}{n}) T$ for $n\in \mathbb{N}$ large enough so that $h_n$ satisfies (\ref{eqn:condh}). It is straightforward to check that $\mu_{h_n ,f} (\lambda) \rightarrow \mu_{h,f} (\lambda)$ as $n\rightarrow +\infty$, and that $\mu_{h_n ,f} (0) <0$ for all $n\in \mathbb{N}$. Thus, for a sufficiently large $n$, we have that $\mu_{h_n ,f} (0) <0$ but $\mu_{h_n ,f} (\lambda) - \lambda^2 >0$, and those two conditions are not equivalent.

\subsection{Boundedness of a sequence of solutions for different speeds}

We first show the following general lemma, which holds without any hypothesis on $\mu_{h,f}$:

\newtheorem{lemma}{Lemma}
\begin{lemma}\label{lm:Tnbounded}
Let $(c_n,T_n,Y_n)$ be a sequence of solutions of $(\ref{eqn:sysfront})$-$(\ref{eqn:condinfty})$ and $(\ref{eqn:neumann})$ such that $0<T_n$ and $0<Y_n <1$ in $\overline{\Omega}$ for each $n \in \mathbb{N}$, and $\sup_n c_n < +\infty$. Then $$\sup_n \|T_n \|_{L^{\infty} (\Omega)} < +\infty.$$
\end{lemma}
\textbf{Proof.} Under those hypotheses, since $c_n \geq c^* $ and $c_n > 0$ for each $n \in \mathbb{N}$ by Theorem \ref{th:thqualit}, we have that the sequence $c_n$ is bounded. Thus, up to extraction of a subsequence, one can assume that $c_n \rightarrow c_\infty \in [\max (c^* ,0) , +\infty )$ as $n \rightarrow +\infty$.

Furthermore, Theorem \ref{th:thqualit} also implies that for each $n \in \mathbb{N}$, the function $T_n$ is globally bounded. Assume now, for the sake of a contradiction, that the sequence $(\|T_n\|_{L^\infty (\Omega)})_{n\in \mathbb{N}}$ is not bounded. Up to extraction of a subsequence, one can assume that $\|T_n\|_{L^\infty (\Omega)}~\rightarrow~+\infty$ as $n\rightarrow +\infty$.

From the boundary conditions (\ref{eqn:condinfty}) and Theorem \ref{th:thqualit}, we know that each pair $T_n$ satisfies $T_n (-\infty ,.) = T_n (+\infty ,.) =0$. Thus, each $T_n$ attains a maximum inside the cylinder $\overline{\Omega}$, and there exists a sequence of points $(x_n ,y_n) \in \overline{\Omega}$ such that
$$T_n (x_n ,y_n) = \max_{\overline{\Omega}} T_n \rightarrow +\infty \mbox{ as } n\rightarrow +\infty .$$
Up to extraction of another subsequence, we may assume that $y_n \rightarrow y_\infty \in \overline{\omega}$ as $n \rightarrow +\infty$.
Define now the normalized shifts
$$U_n (x,y) = \frac{T_n (x+x_n ,y)}{T_n (x_n ,y_n)} .$$
Each function $U_n$ satisfies $0 < U_n \leq 1$ in $\overline{\Omega}$ and is a solution of
$$
\left\{
\begin{array}{rcll}
\displaystyle \Delta U_n + (c_n - u(y)) U_{n,x} + \frac{f(y,T_n (x_n ,y_n)U_n)}{T_n (x_n ,y_n)} Z_n - g_n U_n & = & 0 & \mbox{in } \Omega ,\\
\displaystyle \frac{\partial U_n }{\partial n} & = & 0 & \mbox{on } \partial \Omega ,\\
\end{array}
\right.
$$
where $$Z_n (x,y) = Y_n (x+x_n ,y)$$ is the shifted concentration, and
$$g_n = \frac{h(y, T_n (x+x_n ,y ))}{T_n (x+x_n ,y)} .$$
We already saw in Section \ref{sec:Tbounded} that from the bounds on $h$ in (\ref{eqn:condh}), we have that the sequence $(g_n)_{n\in \mathbb{N}}$ is bounded in $L^\infty (\Omega)$ and thus, up to extraction of a subsequence, we can assume that $g_n$ converges to a function $g$ weakly in $L^{1,*} (\overline{\Omega})$ as $n \rightarrow +\infty$. Furthermore, one can easily check that for all $n \in \mathbb{N}$ and $(x,y)\in \overline{\Omega}$, we have $g_n (x,y) \geq \frac{ \partial h}{\partial T} (y,0)$, whence $g$ is nonnegative and positive on a set of positive measura.
\\

In order to pass to the limit as $n \rightarrow +\infty$, we now claim that
\begin{equation}\label{eqn:lemme1}
\forall K \subset \overline{\Omega} \mbox{ compact, } \max_{(x,y) \in K} Z_n (x ,y) \rightarrow 0 \mbox{ as } n\rightarrow +\infty .
\end{equation}
Indeed, since $Y_n$, $f(y,T_n)/T_n$ and $h(y,T_n)/T_n$ are bounded in $\Omega$ uniformly with respect to $n \in \mathbb{N}$, it follows from Harnack inequality up to the boundary that
\begin{equation}\label{eqn:liminftyTn}
T_n (x+x_n ,y) \rightarrow +\infty \mbox{ as } n\rightarrow +\infty \mbox{ locally uniformly in } (x,y) \in \overline{\Omega} .
\end{equation}
Then, let $K$ be any compact set in $\overline{\Omega}$, and $a \geq 1$ such that $K \subset [-a +1 ,a-1] \times \overline{\omega}$. Define also, for each $n \in \mathbb{N}$:
$$M= \sup_{n \in \mathbb{N} ,y \in \overline{\omega}} | c_n - u(y) | < +\infty ,$$
and
$$m_n = \min_{(x,y) \in [-a ,a] \times \overline{\omega}} f(y,T_n (x+x_n ,y)) \in (0,+\infty ) .$$
From (\ref{eqn:liminftyTn}) and the fact that $f(y,+\infty )= +\infty$ uniformly in $y \in \overline{\omega}$, we have that $m_n \rightarrow +\infty$ as $n \rightarrow +\infty$. Define now, for each $n \in \mathbb{N}$,
$$\lambda_n = \frac{-M + \sqrt{M^2 + 4 \text{Le}^{-1} m_n}}{2 \text{Le}^{-1}} >0$$
the positive solution of
$$\text{Le}^{-1} \lambda_n^2 + M \lambda_n - m_n = 0 .$$
Note that $\lambda_n \rightarrow +\infty$ as $n \rightarrow + \infty$. Lastly, we define
$$\overline{Z}_n (x,y) = e^{-\lambda_n (x+a) } + e^{-\lambda_n (-x+a)} .$$
We now show that $\overline{Z}_n$ is a super-solution for the shifted concentration $Z_n$ in the domain $\Omega_a = (-a ,a) \times \omega$. Both $\overline{Z}_n$ and $Z_n$ satisfy the Neumann boundary conditions on $\partial \Omega$ while
$$Z_n (\pm a,.) \leq 1 \leq \overline{Z}_n (\pm a ,.) \mbox{ in } \overline{\omega} .$$
Inside the domain $\Omega_a$, the function $Z_n$ satisfies
$$0 = \text{Le}^{-1} \Delta Z_n + (c_n - u(y))Z_{n,x} - f(y,T_n (x+x_n ,y))Z_n \leq \text{Le}^{-1}\Delta Z_n + (c_n - u(y)) Z_{n,x} - m_n Z_n ,$$
while $\overline{Z}_n$ satisfies
$$\text{Le}^{-1} \Delta \overline{Z}_n + (c_n - u(y)) \overline{Z}_{n,x} - m_n \overline{Z}_n \leq (\text{Le}^{-1} \lambda_n^2 + M\lambda_n - m_n) \overline{Z}_n = 0 ,$$
owing to the definition of $\lambda_n$. The weak maximum principle then yields
$$ 0 \leq Z_n \leq \overline{Z}_n$$
in $\Omega_a$, for each $n \in \mathbb{N}$. Since $K \subset [ -a +1 , a -1] \times \overline{\omega}$ and $\lambda_n \rightarrow +\infty$ as $n \rightarrow + \infty$, it follows from the definition of $\overline{Z}_n$ that
$$\max_{(x,y) \in K} Z_n (x ,y) \leq \max_{(x,y) \in K} \overline{Z}_n (x,y) \rightarrow 0 \mbox{ as } n\rightarrow +\infty ,$$
and the proof of the claim (\ref{eqn:lemme1}) is now complete.
\\

Lastly, we know that the functions $U_n$ are uniformly bounded (by 1) in $L^\infty (\Omega)$, that
$$0< \frac{f(y,T_n (x_n ,y_n )U_n)}{T_n (x_n ,y_n )} \leq \frac{\partial f}{\partial T}(y,0) U_n \leq \max_{y\in \overline{\omega}} \frac{\partial f}{\partial T}(y,0) \mbox{ in } \Omega ,$$
and that the sequence $(g_n)_{n\in \mathbb{N}}$ is bounded in $L^\infty (\Omega)$. Therefore, we can conclude by standard elliptic estimates that up to extraction of some subsequence, the functions $U_n$ converge as $n \rightarrow +\infty$ in $W_{loc}^{2,p} (\overline{\Omega})$ weak for all $1<p<+\infty$ and strongly in $C_{loc}^1 (\overline{\Omega})$ to a function $U_\infty$ which satisfies
$$
\left\{
\begin{array}{rcll}
\displaystyle \Delta U_\infty + (c_\infty - u(y)) U_{\infty ,x} - g U_\infty & = & 0 & \mbox{in } \Omega , \\
\displaystyle \frac{\partial U_\infty }{\partial n} & = & 0 & \mbox{on } \partial \Omega .\\
\end{array}
\right.
$$
Furthermore, $0 \leq U_\infty \leq 1$ and $U_\infty (0,y_\infty ) =1$. The strong maximum principle and the Hopf lemma then imply that $U_\infty =1$ in $\overline{\Omega}$. This is a contradiction, since $g$ is positive on a set of positive measura. The lemma is now proved. $\Box$

\subsection{Characterization of $Y(-\infty,.)$}\label{sec:lemmaYinfty}

We now show the following lemma, which also holds without any hypothesis on $\mu_{h,f}$:

\begin{lemma}\label{lm:Yinfty}
Let $(c,T,Y)$ be a solution of $(\ref{eqn:sysfront})$ and $(\ref{eqn:neumann})$ such that $0<T$, $0<Y<1$, $T(-\infty ,.) =0$ and $Y(-\infty ,.)$ exists. Then there exists $\beta \geq 0$ such that $\mu_{h,Y_\infty f} (-\beta ) = c\beta + \beta^2$ where $Y_\infty = Y(-\infty ,.)$.
\end{lemma}

\begin{remark}
\upshape Note that, by Theorem \ref{th:thqualit}, any solution of $(\ref{eqn:sysfront})$-$(\ref{eqn:condinfty})$ and $(\ref{eqn:neumann})$ such that $0<T$ and $0<Y<1$ satisfies this lemma.
\end{remark}
\textbf{Proof of Lemma 2.} By Harnack's inequality, we know that $|\nabla T| /T$ is globally bounded. Let
$$\displaystyle \beta = \limsup_{x\rightarrow -\infty} \ \max_{y \in \overline{\omega}} \frac{T_x (x,y)}{T(x,y)},$$
and let us check that $\beta$ satisfies the conclusion of the lemma. First, since $T(-\infty ,.) =0$ and $T>0$, $\beta$ is nonnegative. Let $(x_k ,y_k)_{k \in \mathbb{N}}$ be a sequence of points in $\mathbb{R} \times \overline{\omega}$ such that $x_k \rightarrow -\infty$ and $T_x (x_k ,y_k )/T(x_k ,y_k ) \rightarrow \beta$ as $k \rightarrow +\infty$, and set
$$T_k (x,y) = \frac{T(x_k +x,y)}{T(x_k ,y_k )} .$$
The functions $T_k$ are locally bounded in $\overline{\Omega}$, while the functions $(x,y) \mapsto T(x_k +x,y)$ converge to 0 locally uniformly as $k\rightarrow +\infty$. Therefore, the functions $T_k$ are bounded in all $W_{loc}^{2,p} (\overline{\Omega})$ for all $1\leq p < +\infty$ and converge, up to extraction of a subsequence, to a solution $T_\infty$ of
$$\left\{
\begin{array}{rcll}
\displaystyle \Delta T_\infty + (c-u(y)) T_{\infty ,x} + (\frac{\partial f}{\partial T} (y,0) Y_\infty - \frac{\partial h}{\partial T} (y,0)) T_\infty & = & 0 & \mbox{in } \Omega ,\\
\displaystyle \frac{\partial T_\infty }{\partial n} & = & 0 & \mbox{on } \partial \Omega .\\
\end{array}
\right.
$$
One can also assume that $y_k \rightarrow y_\infty \in \overline{\omega}$. The nonnegative function $T_\infty$ satisfies $T_\infty (0,y_\infty )=1$, whence $T_\infty >0$ in $\overline{\Omega}$ from the strong maximum principle and Hopf lemma. Furthermore, the function
$$z(x,y) = \frac{T_{\infty ,x} (x,y)}{T_\infty (x,y)}$$
satisfies $z \leq \beta$, $z(0,y_\infty ) = \beta$ and
$$\left\{
\begin{array}{rcll}
\displaystyle \Delta z + 2 \frac{\nabla T_\infty}{T_\infty}.\nabla z + (c-u(y)) z_x & = & 0 & \mbox{in } \Omega ,\\
\displaystyle \frac{\partial z }{\partial n} & = & 0 & \mbox{on } \partial \Omega .\\
\end{array}
\right.
$$
The strong maximum principle and Hopf lemma then yield $z= \beta$ in $\overline{\Omega}$. In other words, there exists a positive function $\phi$ in $\overline{\omega}$ such that $T_\infty (x,y) = e^{\beta x} \phi (y)$. The function $\phi$ satisfies
$$\left\{
\begin{array}{rcll}
\displaystyle \Delta \phi +\beta^2 \phi + \beta (c-u(y)) \phi + \frac{\partial f}{\partial T} (y,0) Y_\infty \phi - \frac{\partial h}{\partial T} (y,0) \phi & = & 0 & \mbox{in } \omega ,\\
\displaystyle \frac{\partial \phi }{\partial n} & = & 0 & \mbox{on } \partial \omega .\\
\end{array}
\right.
$$
By uniqueness of the principal eigenvalue for problem (\ref{eqn:principaleigenvalueh}), one can concludes that $\mu_{h,Y_\infty f} (-\beta )= c\beta +\beta^2$. The proof of the lemma is now complete. $\Box$
\\

Let now $\phi_{-\beta }$ be the principal eigenfunction of (\ref{eqn:principaleigenvalueh}) normalized so that $\| \phi_{-\beta}\|_{L^2 (\omega)} =1$. One can easily check that
$$\mu_{h,Y_\infty f} (-\beta) \leq \mu_{h,f} (-\beta) + \int_\omega (1 -Y_\infty ) \frac{\partial f}{\partial T} (y,0) \phi_{-\beta}^2 (y) dy ,$$ and thus
\begin{eqnarray*} Y_\infty \int_\omega \frac{\partial f}{\partial T} (y,0) \phi_{-\beta}^2 (y) dy & \leq & \mu_{h,f} (-\beta) -c\beta -\beta^2+ \int_\omega \frac{\partial f}{\partial T} (y,0) \phi_{-\beta}^2 (y) dy \\
& \leq & \mu_{h,f} (0) + \int_\omega \frac{\partial f}{\partial T} (y,0) \phi_{-\beta}^2 (y) dy , \\
\end{eqnarray*}
where the last inequality follows from the concavity of $\mu_{h,f}$: indeed, for $c \geq c^*$, we have that $c \geq \mu '_{h,f} (0)$, and thus $\mu_{h,f} (\lambda) - (\lambda^2 - c\lambda ) \leq \mu_{h,f} (0)$ for all $\lambda \geq 0$. Note that we already know that $Y_\infty <1$. Nevertheless, when $\mu_{h,f} (0)<0$, we have obtained here a new a priori upper bound on $Y_\infty$, that is
\begin{equation}\label{eqn:a**}
\displaystyle Y_\infty \leq a^* := 1 + \frac{\displaystyle \mu_{h,f} (0 ) }{\displaystyle \int_\omega \frac{\partial f}{\partial T} (y,0) \phi_{-\beta}^2 (y) dy} <1 .\end{equation}
\begin{remark}
\upshape We could also use the fact that $0 \leq c\beta + \beta^2 =\mu_{h,Y_\infty f} (-\beta) $ and $\mu_{h,Y_\infty f} (-\beta) \leq |\omega|^{-1} \int_\omega (\frac{\partial h}{\partial T} (y,0) - Y_\infty \frac{\partial f}{\partial T} (y,0)) dy$ to obtain a more explicit upper bound on $Y_\infty$, that is $$\displaystyle Y_\infty \leq |\omega |^{-1} \ \displaystyle \frac{\displaystyle \int_\omega \frac{\partial h}{\partial T} (y,0) dy}{\displaystyle \int_\omega \frac{\partial f}{\partial T} (y,0) dy} \ .$$
\end{remark}

\subsection{Proof of part (b) of Theorem \ref{th:thexist}}

We now assume that $\sup_{\lambda \in \mathbb{R}} (\mu_{h,f} (\lambda ) - \lambda^2 )< 0 .$ Note first that it immediately follows from elementary geometric considerations that $c^* >0$.
To prove the existence of a non trivial travelling front solution with speed $c^*$, we use an approximation by a sequence of fronts with speeds larger than $c^*$ that we have already constructed.

To do this, let $(c_n)_{n \in \mathbb{N}}$ be a sequence of speeds such that $c_n > c^*$ for all $n$, and such that
$$c_n \rightarrow c^* \mbox{ as } n\rightarrow +\infty .$$
It follows from the results of Section \ref{sec:exist} that for each $n$, there exists a bounded solution $(T_n ,Y_n)$ of (\ref{eqn:sysfront})-(\ref{eqn:condinfty}) and (\ref{eqn:neumann}) with the speed $c = c_n$, such that $T_n >0$ and $0 < Y_n <1$ in $\overline{\Omega}$. According to (\ref{eqn:condinfty}) and Theorem \ref{th:thqualit}, we have
$$T_n (+\infty ,.)=0 \mbox{ and } Y_n (+\infty ,.)=1 ,$$
$$T_n (-\infty ,.)=0 \mbox{ and } Y_n (-\infty ,.)= Y_{n,\infty} \in (0,1) .$$
It also follows from Lemma \ref{lm:Tnbounded} that there exists a constant $M>0$ such that
\begin{equation}\label{eqn:Tbounded}
\forall n \in \mathbb{N} \mbox{, } \forall (x,y) \in \overline{\Omega} \mbox{, } 0 < T_n (x,y) \leq M .
\end{equation}
As we have mentioned, our strategy is to pass to the limit as $n \rightarrow +\infty$, in order to get a solution of (\ref{eqn:sysfront})-(\ref{eqn:condinfty}) and (\ref{eqn:neumann}) with the speed $c = c^*$. Any shift of the travelling wave $(T_n , Y_n)$ in the variable $x$ along the cylinder is, of course, also a travelling wave, and the main technical difficulty here is to shift suitably the functions $(T_n ,Y_n )$ so that the limit pair is non-trivial and satisfies the correct limiting conditions at infinity. For that we have to identify a region where both $T_n$ and $Y_n$ are uniformly not very flat.

\subsubsection*{Locating the interface}

Let $a^*$ defined in (\ref{eqn:a**}). For each $a \in (a^*,1)$, and $n \in \mathbb{N}$, we define
$$x_n^a = \min \{ x \in \mathbb{R} , Y_n \geq a \mbox{ in } [x,+\infty ) \times \overline{\omega} \} .$$
Since the functions $Y_n$ are continuous in $\overline{\Omega}$, satisfy $Y_n (+\infty ,.) =1$ and $Y_n (-\infty ,.) \leq a^* $ by (\ref{eqn:a**}), the real numbers $x_n^a$ are well-defined. Moreover, $x_n^a$ is nondecreasing in $a \in (a^*,1)$ for each fixed $n$. Observe that, also,
$$\left\{
\begin{array}{l}
\displaystyle Y_n \geq a \mbox{ in } [x_n^a , +\infty ) \times \overline{\omega} ,\\
\displaystyle \min_{\overline{\omega}} Y_n (x_n^a ,.) = a .\\
\end{array}
\right.
$$
Since $Y_n (+\infty ,.) =1$, we have
$$\| \nabla Y_n \|_{L^\infty ([x_n^a ,+ \infty ) \times \overline{\omega })} >0 .$$
Furthermore, since $| \nabla Y_n (x,y) | \rightarrow 0$ as $x \rightarrow +\infty$ uniformly in $y \in \overline{\omega}$, the points
$$\tilde{x}_n^a = \min \{x \in [x_n^a , +\infty ) \mbox{, } \exists y \in \overline{\omega} \mbox{, } | \nabla Y_n (x,y) | = \| \nabla Y_n \|_{L^\infty ([x_n^a ,+ \infty ) \times \overline{\omega })} \}$$
are well-defined.

We now introduce the following lemma, that shows that to the right of $x_n^a$, there are regions where $Y_n$ are uniformly non too flat.

\begin{lemma}\label{lm:nablaYn}
For all $a \in (a^*,1)$, we have
$$\inf_{n} \| \nabla Y_n \|_{L^\infty ([x_n^a ,+ \infty ) \times \overline{\omega })} >0 .$$
\end{lemma}
The proof of this lemma is postponed until the end of the section.

\subsubsection*{Normalization of $(T_n ,Y_n )$ and passage to the limit}

Let us now complete the proof of the existence of a non-trivial bounded solution $(T,Y)$ of (\ref{eqn:sysfront})-(\ref{eqn:condinfty}) and (\ref{eqn:neumann}) with the speed $c=c^*$.

Choose now any $a \in (a^*,1)$ and let $\tilde{y}_n^a$ be a sequence of points in $\overline{\omega}$ such that
$$|\nabla Y_n (\tilde{x}_n^a ,\tilde{y}_n^a ) | = \| \nabla Y_n \|_{L^\infty ([x_n^a ,+ \infty ) \times \overline{\omega })}$$
for all $n \in \mathbb{N}$. Up to extraction of a subsequence, one can assume that the sequence $\tilde{y}_n^a$ converges to a point $\tilde{y}^a \in \overline{\omega}$. Lemma \ref{lm:nablaYn} implies that
\begin{equation}\label{eqn:infYn}
\inf_n | \nabla Y_n (\tilde{x}_n^a ,\tilde{y}_n^a ) | >0 .
\end{equation}
For each $n$ and $(x,y) \in \overline{\Omega}$, define the shifted functions
$$T_n^a (x,y) = T_n (x+\tilde{x}_n^a ,y) ,$$
$$Y_n^a (x,y) = Y_n (x+\tilde{x}_n^a ,y) .$$
Recall that both $T_n$ and $Y_n$ are uniformly bounded in $\overline{\Omega}$, independently of $n$ (that is (\ref{eqn:Tbounded})). By standard elliptic estimates up to the boundary, these functions, as well as the shifts $T_n^a$ and $Y_n^a$ are also bounded in $C^{2,\alpha} (\overline{\Omega})$, uniformly in $n$. Up to extraction of a subsequence, one can assume that the sequence $( T_n^a , Y_n^a )$ converges to a function $(T^a ,Y^a )$ in $C_{loc}^2 (\overline{\Omega})$ as $n\rightarrow +\infty$. Passing to the limit, we conclude that the pair $(T^a ,Y^a )$ satisfies
\begin{equation}\label{eqn:sysfrontc*}
\left\{
\begin{array}{rcll}
\Delta T^a + (c^* -u(y)) T_x^a + f(y,T^a )Y^a - h(y,T^a ) & = & 0 & \mbox{in } \Omega ,\\
\text{Le}^{-1} \Delta Y^a + (c^* - u(y)) Y_x^a - f(y,T^a )Y^a & = & 0 & \mbox{in } \Omega ,\\
\end{array}
\right.
\end{equation}
with the Neumann boundary conditions on $\partial \Omega$
$$\frac{\partial T^a }{\partial n} = \frac{\partial Y^a }{\partial n} =0 ,$$
and they obey the uniform bounds $0 \leq T^a \leq M$ and $0 \leq Y^a \leq 1$ in $\overline{\Omega}$. Furthermore, (\ref{eqn:infYn}) implies that
\begin{equation}\label{eqn:Ynontrivial}
| \nabla Y^a (0, \tilde{y}^a ) | >0 .
\end{equation}
Thus, $Y^a$ is not a constant. By the strong maximum principle and Hopf lemma, we can conlude that $0< Y^a < 1$ in $\overline{\Omega}$, and $Y^a$ is non-trivial.

Let us now check that $T^a >0$. Otherwise, if $T^a$ vanishes somewhere in $\overline{\Omega}$, then it is identically equal to 0 by the strong maximum principle and Hopf Lemma. In that case, the function $Y^a$ would satisfy
\begin{equation}\label{eqn:ifTvanish}
\left\{
\begin{array}{rcll}
\displaystyle \text{Le}^{-1} \Delta Y^a + (c^* -u(y)) Y_x^a & = & 0 & \mbox{in } \Omega ,\\
\displaystyle \frac{\partial Y^a}{\partial n} & = & 0 & \mbox{on } \partial \Omega . \\
\end{array}
\right.
\end{equation}
We apply now the same method as in Sections \ref{sec:Tbounded} and \ref{sec:behavleft}. If we multiply (\ref{eqn:ifTvanish}) by $Y^a$, integrate over a finite cylinder $(-A,A) \times \omega$ and pass to the limit as $A \rightarrow + \infty$, we would obtain that the integral
$$\int_\Omega | \nabla Y^a |^2$$
is finite. Then, for a sequence $A_n \rightarrow +\infty$, the shifted functions $Y^a (\pm A_n +x ,y)$ would converge in $C_{loc}^2 (\overline{\Omega})$ to two constants $Y_{\pm }^a \in [0,1]$. Integrating (\ref{eqn:ifTvanish}) over the cylinder $(-A_n ,A_n ) \times \omega$ and passing to the limit as $n \rightarrow +\infty $ then yield that
$$c^* (Y_+^a - Y_-^a ) = 0,$$
that is $$Y_+^a = Y_-^a .$$
Finally, once again, multiplying (\ref{eqn:ifTvanish}) by $Y^a$, integrating over the cylinder $(-A_n,A_n) \times \omega$ and passing to the limit as $n \rightarrow + \infty$ imply that
$$\int_\Omega | \nabla Y^a |^2 = 0 ,$$
which contradicts (\ref{eqn:Ynontrivial}). We conclude that $T^a >0$ in $\overline{\Omega}$.

\subsubsection*{The limits at infinity}

It only remains to show that $T^a$ and $Y^a$ attain the correct limits at infinity. As before, we can show that the integrals
$$\int_\Omega f(y,T^a)Y^a \mbox{, } \int_\Omega \frac{\partial h}{\partial T}(y,0) T^a \mbox{, } \int_\Omega | \nabla T^a |^2  \mbox{ and } \int_\Omega | \nabla Y^a |^2$$
converge. Therefore, for any sequence $A_n \rightarrow +\infty$, there exists a subsequence such that the functions $T^a (x\pm A_n ,y)$ and $Y^a (x \pm A_n ,y)$ converge in $C_{loc}^2 (\overline{\Omega})$ as $n \rightarrow +\infty$ to some nonnegative constants $T_{\pm}$ and $Y_{\pm}$. We then have that
$$\int_\omega \frac{\partial h}{\partial T}(y,0) T_\pm dy= 0,$$ thus $T_{\pm} =0$ independent of the sequence $(A_n)_{n\in\mathbb{N}}$, and $T^a(x,y) \rightarrow 0$ as $x \rightarrow \pm \infty$ uniformly in $y \in \overline{\omega}$. By standard elliptic estimates, we also have that $T_x^a (-\infty ,.) =0$.

Furthermore, let $(A_n)_{n\in\mathbb{N }}$ and $(B_n)_{n \in\mathbb{N}}$ be two sequences which converge to $+\infty$ as $n\rightarrow +\infty$, and such that $Y_+ := \lim_n Y^a (x + A_n ,y)$ and $Y_- :=\lim_n Y^a (x - B_n ,y)$ are well-defined. By integrating the equation (\ref{eqn:sysfrontc*}) satisfied by $Y^a$ over the domain $(-B_n , A_n ) \times \omega$ and passing to the limit $n \rightarrow +\infty$, we obtain
$$\int_\Omega f(y,T^a ) Y^a = c^* (Y_+ - Y_- ) |\omega | .$$
If we fix the sequence $(B_n)_{n \in \mathbb{N}}$, we see that $Y_+$ does not depend on the choice of the sequence $(A_n)_{n \in \mathbb{N}}$. Symmetrically, we have that $Y_-$ does not depend on the choice of $(B_n)_{n \in \mathbb{N}}$. By standard elliptic estimates, we also deduce that $Y_x^a (\pm \infty ,.) = 0$.

We will now show that $Y_+ =1$. We first claim that the sequence $z_n^a~=~\tilde{x}_n^a -x_n^a$ is bounded. Otherwise, up to extraction of another subsequence, we would have $z_n^a \rightarrow +\infty$ as $n \rightarrow +\infty$. Thus, for each $(x,y) \in \overline{\Omega}$, we would have $x + \tilde{x_n^a} \geq x_n^a$ for sufficiently large $n$, and so
$$Y_n^a (x,y) = Y_n (x+\tilde{x}_n^a ,y) \geq a$$
for sufficiently large $n$. This would imply that $Y^a (x,y) \geq a > a^*$ in $\overline{\Omega}$. But $Y^a (-\infty ,.) \leq a^*$ from the calculations in Section \ref{sec:lemmaYinfty}, which leads to a contradiction.

Let now $b$ be any real number in $(a,1)$. As in the previous argument, the shifted functions
$$T_n^b (x,y) =T_n (x+ \tilde{x}_n^b ,y) \mbox{,  } Y_n^b (x,y) = Y_n(x+ \tilde{x}_n^b ,y)$$
converge in $C_{loc}^2 (\overline{\Omega})$ as $n \rightarrow +\infty$, up to extraction of a subsequence, to a pair $(T^b ,Y^b )$ of solutions of (\ref{eqn:sysfrontc*}) with $b$ instead of $a$, such that $Y^b (-\infty ,.) \leq a^*$ from the calculations in Section \ref{sec:lemmaYinfty}. We claim that the sequence $(x_n^b -x_n^a)_{n\in \mathbb{N}}$ is bounded. Indeed, as we know that the sequence $z_n^b =\tilde{x}_n^b - x_n^b$ is bounded, if the sequence $(x_n^b -x_n^a)$ is unbounded, then the sequence of nonnegative numbers $(\tilde{x}_n^b -x_n^a)$ would be unbounded, which would imply that $Y^b~(-\infty ,.)~\geq~a~> a^*$. This is a contradiction.

Therefore, the sequence $(x_n^b - \tilde{x_n^a})$ is also bounded, and there exists $A_a^b \geq 0$ (which depends on $a$ and $b$ but not on $n$) such that $x_n^b - \tilde{x_n^a} \leq A_a^b$ for all $n$. However, for each $(x,y) \in [A_a^b ,+\infty ) \times \overline{\omega}$, we have then $x + \tilde{x}_n^a \geq x_n^b$ and thus
$$Y_n^a (x,y) = Y_n (x+\tilde{x}_n^a ,y ) \geq b$$
for all $n$. As a consequence, we have that $Y^a (+ \infty ,.) \geq b$.

Since $b$ was arbitrarily chosen in $(a,1)$ and since $Y^a (+\infty ,.) \leq 1$, we conclude that $Y^a (+\infty ,.) =1$. This completes the proof of Theorem \ref{th:thexist}. $\Box$

\subsection{Proof of lemma \ref{lm:nablaYn}}

Assume by contradiction that the conclusion of lemma \ref{lm:nablaYn} does not hold for a real number $a \in (a^* ,1)$. As $\| \nabla Y_n \|_{L^\infty ([x_n^a ,+\infty ) \times \overline{\omega})}$ is positive for each $n \in \mathbb{N}$, up to extraction of a subsequence, one can then assume without loss of generality that
\begin{equation}\label{eqn:contradlemma}
\| \nabla Y_n \|_{L^\infty ([x_n^a ,+\infty ) \times \overline{\omega})} \rightarrow 0 \mbox{ as } n\rightarrow +\infty .
\end{equation}

\subsubsection*{Temperature is small on the right}

We first claim that in this case, the "temperature interface" is located far to the left of the "concentration interface", that is, we have
\begin{equation}\label{eqn:Tsmall}
\| T_n \|_{L^\infty ([x_n^a ,+\infty ) \times \overline{\omega})} \rightarrow 0 \mbox{ as } n\rightarrow +\infty.
\end{equation}
Indeed, assume now that (\ref{eqn:contradlemma}) holds and (\ref{eqn:Tsmall}) does not. Then there exist $\varepsilon >0$ and a sequence $(x_n ,y_n)_{n \in\mathbb{N}}$ in $\overline{\Omega}$ such that
$$x_n \geq x_n^a \mbox{ and } T_n (x_n ,y_n )\geq \varepsilon \mbox{ for all } n \in \mathbb{N} .$$
Up to extraction of a subsequence, we can assume that $y_n \rightarrow y_\infty \in \overline{\omega}$ as $n \rightarrow +\infty$. The standard elliptic estimates imply that the sequence of shifted functions $T_n (x+x_n ,y)$ and $Y_n (x+x_n ,y)$ converge in $C_{loc}^2 (\overline{\Omega})$, up to extraction of some subsequence, to a pair $(T,Y)$ of solutions of (\ref{eqn:sysfrontc*}). Furthermore, $T$ and $Y$ satisfy
$$0 \leq Y \leq 1 \mbox{, } 0 \leq T \leq M \mbox{ in } \overline{\Omega} ,$$
$$Y \geq a >0 \mbox{, } | \nabla Y | = 0 \mbox{ in } [0,+\infty ) \times \overline{\omega} ,$$
and $T(0, y_\infty ) \geq \varepsilon$. The strong maximum principle and Hopf lemma then imply that $T>0$ and $Y>0$ in $\overline{\Omega}$. This is a contradiction because $Y$ is a constant in $[0 ,+\infty ) \times \overline{\omega}$ and thus has to satisfy $f(y,T)Y =0$ in the same domain.

\subsubsection*{Temperature decays exponentially on the right}

We then claim that under assumptions (\ref{eqn:contradlemma}) and hence (\ref{eqn:Tsmall}), $T_n$ decays exponentially uniformly to the right of $x_n^a$, that is: there exist a positive number $\lambda >0$, an integer $N$ and $A\geq 0$ so that for all $n \geq N$ and all $(x,y) \in [x_n^a + A,+\infty ) \times \overline{\omega}$ we have
\begin{equation}\label{eqn:Tdecayexp}
\frac{T_{n,x} (x,y)}{T_n (x,y)} \leq -\lambda .
\end{equation}
As $T_n >0$, while $Y_n$, $f(y,T_n)/T_n$ and $h(y,T_n)/T_n$ are bounded independently of $n$, and $(T_n,Y_n)$ satisfy (\ref{eqn:sysfront}) with the speeds $c_n$ which are uniformly bounded (since $c_n \rightarrow c^*$ as $n\rightarrow +\infty$), it follows from standard elliptic estimates and the Harnack inequality that the functions $| \nabla T_n |/T_n$ are uniformly bounded in $\Omega$. Assume now that the claim (\ref{eqn:Tdecayexp}) does not hold. Then, after extraction of a subsequence, there exists a sequence of points $(x_n ,y_n) \in [x_n^a ,+\infty) \times \overline{\omega}$ such that
\begin{equation}\label{eqn:eqq1}
\lim_{n\rightarrow +\infty} (x_n - x_n^a ) = +\infty
\end{equation}
and
\begin{equation}\label{eqn:eqq2}
\liminf_{n\rightarrow +\infty} \frac{T_{n,x} (x_n ,y_n )}{T_n (x_n ,y_n )} \geq 0 .
\end{equation}
Set the normalized and shifted temperature
$$U_n (x,y) = \frac{T_n (x+x_n ,y)}{T_n (x_n ,y_n)}$$
for all $n$ and $(x,y) \in \overline{\Omega}$. Up to extraction of another subsequence, one can assume that $y_n \rightarrow y_\infty \in \overline{\omega}$ as $n \rightarrow +\infty$. The functions $U_n$ satisfy
$$
\left\{
\begin{array}{rcll}
\displaystyle \Delta U_n + (c_n -u(y)) U_{n,x} + \frac{f(y,T_n (x_n ,y_n ) U_n)}{T_n (x_n ,y_n )} Z_n - \frac{h(y,T_n (x_n ,y_n ) U_n)}{T_n (x_n ,y_n )} & = & 0 & \mbox{in } \Omega ,\\
\displaystyle \frac{\partial U_n }{\partial n} & = & 0 & \mbox{on } \partial \Omega ,\\
\end{array}
\right.
$$
where $Z_n (x,y) = Y_n (x+ x_n ,y)$ is the shifted concentration. The sequence $U_n$ is bounded in $L_{loc}^\infty (\overline{\Omega})$ and in $W_{loc}^{2,p} (\overline{\Omega})$ (for all $1 \leq p < +\infty$) while $T_n (x_n ,y_n) \rightarrow 0$ as $n \rightarrow +\infty$, as can be seen from (\ref{eqn:Tsmall}) because $x_n \geq x_n^a$. On the other hand, the sequence of functions $Z_n$ are globally bounded in $C^{2,\alpha} (\overline{\Omega})$. Hence, up to extraction of a subsequence, the functions $Z_n$ converge to a function $Z$ in $C_{loc}^2 (\overline{\Omega})$ as $n\rightarrow +\infty$. But (\ref{eqn:contradlemma}) and (\ref{eqn:eqq1}) imply that $Z$ is a constant. Furthermore, since $a \leq Y_n \leq 1$ in $[x_n^a ,+\infty )\times \overline{\omega}$, the constant $Z$ is such that $$0 <a \leq Z \leq 1 .$$
As a consequence, up to extraction of another subsequence, the positive functions $U_n$ converge in all $W_{loc}^{2,p} (\overline{\Omega})$ weak (for $1 <p<+\infty$) to a classical nonnegative solution $U$ of
$$
\left\{
\begin{array}{rcll}
\displaystyle \Delta U + (c^* -u(y)) U_{x} + \frac{\partial f}{\partial T} (y,0))Z U - \frac{\partial h}{\partial T} (y,0) U & = & 0 & \mbox{in } \Omega ,\\
\displaystyle \frac{\partial U }{\partial n} & = & 0 & \mbox{on } \partial \Omega .\\
\end{array}
\right.
$$
Furthermore, we have that $U(0, y_\infty ) =1$ while (\ref{eqn:eqq2}) implies
$$\frac{U_x (0 ,y_\infty )}{U (0, y_\infty )} \geq 0 .$$
It follows from the strong maximum principle and the Hopf lemma that $U>0$ in $\overline{\Omega}$, and it follows from standard elliptic estimates and the Harnack inequality that the function $|\nabla U |/U$ is bounded in $\Omega$. Let $(x'_n ,y'_n)_{n\in \mathbb{N}}$ be a sequence of points in $\overline{\Omega}$ such that
$$\frac{U_x (x'_n ,y'_n )}{U(x'_n ,y'_n)} \rightarrow \sup_{\overline{\Omega}} \frac{U_x }{U} =: \overline{M} \geq 0 \mbox{ as } n\rightarrow +\infty .$$
Up to extraction of a subsequence, one can assume that $y'_n \rightarrow y'_\infty$ as $n\rightarrow +\infty$. Next, with the same arguments as above, the functions
$$V_n (x,y) = \frac{U(x+x'_n ,y)}{U(x'_n ,y'_n )}$$
are bounded in $C_{loc}^{2,\alpha} (\overline{\Omega})$ independently of $n$ and converge in $C_{loc}^2 (\overline{\Omega})$, up to extraction of some subsequence, to a nonnegative function $V$ solving the same linear equation as $U$. Furthermore, we have that $V (0, y'_\infty ) =1$. Therefore, by the strong maximum principle and the Hopf lemma, $V$ is positive in $\overline{\Omega}$. Moreover, we have
$$\frac{V_x}{V} \leq \overline{M} \mbox{ in } \overline{\Omega} \mbox{ and } \frac{V_x (0, y'_\infty)}{V(0,y'_\infty)}=\overline{M} .$$
However, one can easily check that the function $V_x /V$ satisfies a linear elliptic equation in $\overline{\Omega}$ without the zeroth-order term, together with the Neumann boundary condition on $\partial \Omega$. Since $V_x /V$ attains its maximum at the point $(0 , y'_\infty)$, the maximum principle implies that $V_x /V$ is identically equal to $\overline{M}$ in $\overline{\Omega}$. In other words, there exists a positive function $\phi (y)$ such that $V(x,y) = e^{\overline{M}x} \phi (y)$ in $\overline{\Omega}$. It follows that the function $\phi$ satisfies
$$
\left\{
\begin{array}{rcll}
\displaystyle \Delta \phi + \overline{M}^2 \phi +\overline{M} (c^* - u(y)) \phi + \frac{\partial f}{\partial T} (y,0)Z \phi - \frac{\partial h}{\partial T} (y,0) \phi & = & 0 & \mbox{in } \overline{\omega} ,\\
\displaystyle \frac{\partial \phi}{\partial n} & = & 0 & \mbox{on } \partial \omega .\\
\end{array}
\right.
$$
By uniqueness of the principal eigenvalue for (\ref{eqn:principaleigenvalueh}), we conclude that
$$\mu_{h,Zf} (-\overline{M}) = c^* \overline{M} + \overline{M}^2 .$$
Recall that $\overline{M} \geq 0$ and $c^* >0$. Hence, as in section \ref{sec:lemmaYinfty}, it implies that $Z \leq a^*$. But we saw that $Z \geq a > a^*$. One has then reached a contradiction which shows that (\ref{eqn:Tdecayexp}) must hold.

\subsubsection*{A sub-solution for $Y_n$}

We have just shown that for all $n \geq N$ and $(x,y) \in [x_n^a + A,+\infty ) \times \overline{\omega}$, we have
$$0< T_n (x,y) \leq T_n (x_n^a +A ,y) e^{-\lambda (x-x_n^a -A)} \leq M e^{-\lambda (x-x_n^a -A)} .$$
The last inequality above follows from (\ref{eqn:Tbounded}). On the other hand, for all $x \in [x_n^a , x_n^a +A]$ we have that $e^{-\lambda (x-x_n^a -A)} \geq 1$. Then the above inequality holds in the whole half-strip $x\geq x_n^a$:
\begin{equation}\label{eqn:Tdecayexp2}
\forall n \geq N \mbox{, } \forall (x,y) \in [x_n^a ,+\infty ) \times \overline{\omega} \mbox{, } 0 <T_n (x,y) \leq M e^{-\lambda (x-x_n^a -A)} .
\end{equation}
We apply the same strategy as in Section \ref{sec:subsolY}: we use the above exponential bound for temperature to create a sub-solution for $Y_n$. First, since $\nu (0) = \nu '(0) = 0 < c^*$, one can choose $\beta >0$ small enough so that
\begin{equation}\label{eqn:betabis}
\left\{
\begin{array}{l}
0<\beta < \lambda ,\\
\nu (\beta \text{Le} )-\beta^2 + c^* \beta \text{Le} > 0 ,\\
\end{array}
\right.
\end{equation}
and $\gamma >0$ large enough so that
\begin{equation}\label{eqn:gammabis}
\left\{
\begin{array}{l}
\displaystyle \gamma \times \min_{\overline{\omega}} \psi _{\beta \text{Le}} \geq 1 , \\
\displaystyle \gamma \text{Le}^{-1} (\nu (\beta \text{Le}) -\beta^2 + c^* \beta \text{Le} ) \times \min_{\overline{\omega}} \psi _{\beta \text{Le}} \geq \max_{y \in \overline{\omega}} \frac{\partial f}{\partial T} (y,0) M e^{\lambda A} ,
\end{array}
\right.
\end{equation}
where $\psi _{\beta \text{Le}}$ denotes the positive principal eigenfunction of (\ref{eqn:principaleigenvalue}) with parameter $\beta \text{Le}$. For each $n \geq N$, we define
$$\underline{Y}_n (x,y) =  \max (0,1-\gamma \psi_{\beta \text{Le}} (y) e^{-\beta (x-x_n^a )})$$
for all $(x,y) \in \overline{\Omega}$. Each function $\underline{Y}_n$ satisfies the Neumann boundary conditions on $\partial \Omega$, while $0 \leq \underline{Y}_n \leq 1$ and $\underline{Y}_n (+\infty ,.) = Y_n (+\infty ,.) =1$ uniformly in $\overline{\omega}$. In addition, it follows from (\ref{eqn:gammabis}) that
$$\underline{Y}_n = 0 \mbox{ in } (-\infty ,x_n^a] \times \overline{\omega} .$$
Therefore, in the region where $\underline{Y}_n (x,y) >0$, we have $x >x_n^a$ and thus there $\underline{Y}_n$ satisfies
$$
\begin{array}{l}
\displaystyle \text{Le}^{-1} \Delta \underline{Y}_n + (c_n -u(y)) \underline{Y}_{n,x} - f(y,T_n)\underline{Y}_n \\
\displaystyle \geq \gamma \text{Le}^{-1} (\nu (\beta \text{Le}) -\beta^2 + c_n \beta \text{Le} ) \psi_{\beta \text{Le}} (y) e^{-\beta (x-x_n^a )}
- \frac{\partial f}{\partial T} (y,0) M e^{-\lambda (x-x_n^a -A)} \\
\displaystyle \geq \gamma \text{Le}^{-1} (\nu (\beta \text{Le}) -\beta^2 + c^* \beta \text{Le} ) \psi_{\beta \text{Le}} (y) e^{-\beta (x-x_n^a )}
- \frac{\partial f}{\partial T} (y,0) M e^{\lambda A} e^{-\beta (x-x_n^a )} \geq 0 ,
\end{array}
$$
because $f$ of the KPP-type, $c_n > c^*$ and from (\ref{eqn:Tdecayexp2})-(\ref{eqn:gammabis}). As $f(y,T_n) \geq 0$, it then follows from the weak maximum principle that
$$\forall n \geq N \mbox{, } \forall (x,y) \in [x_n^a ,+\infty ) \times \overline{\omega} \mbox{,  } Y_n (x,y) \geq \underline{Y}_n (x,y) \geq 1 - \gamma \psi _{\beta \text{Le}} (y) e^{-\beta (x-x_n^a )} .$$
In particular, there exists $L_0 >0$ independent of $n$ so that we have $Y_n (x_n^a + L_0 ,y) \geq (1+a)/2$ for all $y \in \overline{\omega}$. However, since $\min_{\overline{\omega}} Y_n (x_n^a ,y) = a <1$ for all $n$, we finally reach a contradiction to our assumption (\ref{eqn:contradlemma}). This completes the proof of Lemma \ref{lm:nablaYn}. $\Box$

\section{Criteria for flame extinction, blow-off or propagation}\label{sec:cauchy}

This section will deal with the proof of Theorem \ref{th:cauchy} and will be divided in two parts. The first part will treat of both flame extinction and blow-off, which rely on the same method, that is the search for a suitable supersolution for temperature. The case of flame propagation will be treated separately and will use the same method as in Section \ref{sec:exist} to construct not only a supersolution but also a sub-solution.

\subsection{Flame extinction and blow-off}

Let $(T,Y)$ be the solution of the Cauchy problem defined by (\ref{eqn:sys})-(\ref{eqn:neumann}) with an initial profile $(T_0 ,Y_0)$ verifying (\ref{eqn:iniprofile}), and let $\lambda >0$, $C>0$ such that
$$T_0 (x,y) \leq C e^{-\lambda x} \mbox{ in } \mathbb{R}^+ \times \overline{\omega} .$$
Moreover, we have that $T_0$ is bounded. Therefore, by increasing $C$, we can assume without loss of generality that we also have $T_0 \leq C$ in the entire domain $\Omega$.

We then observe that $0 \leq T(t,x,y)$ and $0 \leq Y(t,x,y) \leq 1$ for all $t \geq 0$ and $(x, y ) \in \overline{\Omega}$ as follows from the maximum principle. Therefore, it is straightforward to check that $0 \leq T(t,x,y) \leq \Phi (t,x,y )$ for all $t \geq 0$ and $(x, y ) \in \overline{\Omega}$, where $\Phi$ is any solution of
\begin{equation}\label{eqn:supsolcauchy1}
\left\{
\begin{array}{rl}
\Phi_t + u(y) \Phi_x \geq \Delta \Phi + \frac{\partial f}{\partial T}(y,0) \Phi - \frac{\partial h}{\partial T}(y,0) \Phi & \mbox{for } t\geq 0 \mbox{ and } (x,y) \in \Omega ,\\
\frac{\partial \Phi }{\partial n} = 0 & \mbox{for } t\geq 0 \mbox{ and } (x,y) \in \partial \Omega ,\\
\end{array}
\right.
\end{equation}
provided that $T_0 (x,y) \leq \Phi (0,x,y)$ for all $(x, y ) \in \overline{\Omega}$.

If $\mu_{h,f} (0) > 0$, we choose the function $$\Phi (t,y) = C e^{-\mu_{h,f} (0)t} \phi_0 (y),$$ where $\phi_0$ is the positive eigenfunction of (\ref{eqn:principaleigenvalueh}) with parameter 0, normalized so that $\min_{\overline{\omega}} \phi_0~=~1$. Such a $\Phi$ indeed satisfies (\ref{eqn:supsolcauchy1}), and we have that $T_0 \leq C \leq \Phi(0,y)$, which proves part (a) of Theorem \ref{th:cauchy}. Notice that here, we only used the fact that $T_0$ is bounded.

We now assume that there exists $\eta \in (0, \lambda ]$ such that $\mu_{h,f} (\eta) - \eta^2 > 0$. Let us look for the supersolution $\Phi$ in the form $$\Phi (t,x,y) = C e^{-\eta (x+\gamma t)} \phi_\eta (y) ,$$
where $\gamma >0$ is to be determined. Note that the fact that $\eta \leq \lambda$ and the above bounds on $T_0$ guarantees that $T_0 (x,y) \leq \Phi (0,x,y)$ for all $(x, y ) \in \overline{\Omega}$ regardless of the choice of $\gamma$. Insert now the expression of $\Phi$ in (\ref{eqn:supsolcauchy1}) and obtain that we need
$$-\eta \gamma \phi_\eta - \eta u(y) \phi_\eta \geq \eta^2 \phi_\eta + \Delta_y \phi_\eta + \frac{\partial f}{\partial T}(y,0) \phi_\eta - \frac{\partial h}{\partial T} (y,0) \phi_\eta \ \text{ in } \overline{\omega}.$$
This is true if and only if
$$\eta \gamma \leq \mu_{h,f} (\eta) -\eta^2 .$$
Since the right-hand side is positive, this inequality holds for some small $\gamma >0$, which concludes the proof of part (b) of Theorem \ref{th:cauchy}.

\subsection{Propagation}

Let $(T,Y)$ be the solution of the Cauchy problem defined by (\ref{eqn:sys})-(\ref{eqn:neumann}) with an initial profile $(T_0 ,Y_0)$ verifying (\ref{eqn:iniprofile}), that is there exists $\lambda$, $\lambda ' >0$, and $C_1$, $C_2$, $C_3  >0$ such that
\begin{equation}\label{eqn:iniprofilb}
\begin{array}{l}
0 \leq T_0 \mbox{, } T_0 \mbox{ is bounded, } 0\leq Y_0 \leq 1 , \\
C_1 e^{-\lambda x} \leq T_0 (x,y) \leq C_2 e^{-\lambda x} \mbox{ in } \mathbb{R}^+ \times \overline{\omega},\\
1 - Y_0 (x,y) \leq C_3 e^{-\lambda ' x}  \mbox{ in } \mathbb{R}^+ \times \overline{\omega}.
\end{array}
\end{equation}
We assume that $\mu_{h,f} (0) < 0$, thus $c^*$ and $\lambda^*$ are well defined. We also assume that $k (\lambda) = \lambda^2 - \mu_{h,f} (\lambda) >0$, which implies that $c :=k(\lambda)/\lambda >0$. Lastly, it follows from the hypothesis $\lambda < \lambda^*$ that $c > c^*$ (recall that $k(s)=cs$ has no positive solution for $c < c^*$ and only one for $c= c^*$, which is $\lambda^*$), and that $\lambda$ is the smallest positive root of $k(s)=cs$ (recall that $k(s)=cs$ has two positive solutions $\lambda_1$, $\lambda_2$ for $c > c^*$, with $\lambda_1 < \lambda^* < \lambda_2$).

As before, the maximum principle implies that $0 \leq T(t,x,y) $ and $0 \leq Y(t,x,y) \leq 1$ for all $t \geq 0$ and $(x,y) \in \overline{\Omega}$. We now proceed as in Section \ref{sec:exist} in order to construct sub- and super-solutions for $T$ that both move to speed $c$.

First, we define the function
$$\overline{T} (t,x,y) = C e^{-\lambda (x-ct)} \phi_\lambda (y) >0,$$
where $C>0$ to be determined, and $\phi_\lambda$ is the positive principal eigenfunction of (\ref{eqn:principaleigenvalueh}) with parameter $\lambda$, normalized so that $\| \phi_\lambda \|_{L^\infty (\omega)} =1$. The function $\overline{T}$ satisfies the Neumann boundary conditions on $\partial \Omega$ and
$$\overline{T}_t + u(y) \overline{T}_x = \Delta \overline{T} + \frac{\partial f}{\partial T}(y,0) \overline{T} - \frac{\partial h}{\partial T} (y,0) \overline{T} \geq \Delta \overline{T} + f(y,\overline{T})Y - h(y,\overline{T}) ,$$
for all $t \geq 0$ and $(x,y)\in \Omega$. Furthermore, it follows from (\ref{eqn:iniprofilb}) that we can choose $C$ large enough so that $T_0 (x,y) \leq \overline{T} (0,x,y)$ for all $(x,y) \in \overline{\Omega}$. Therefore, we have $$T(t,x,y) \leq \overline{T} (t,x,y)$$ for all $t\geq 0$ and $(x,y) \in \overline{\Omega}$. This implies in particular that for $(x,y) \in \overline{\Omega}$ and $c' >c$, we have
$$T(t,x+c' t,y) \leq C e^{-\lambda (x+c' t -ct)} \rightarrow 0 \mbox{ as } t\rightarrow +\infty .$$
It now remains to find $x_0 \in \mathbb{R}$ and $\alpha (x_0, y) > 0$ such that $T(t,x_0 + ct ,y) \geq \alpha (x_0, y)$ for all $t \geq 1$ and $y \in \overline{\omega}$. To do this, we search for a suitable sub-solution, as announced above. Let first, as in Section \ref{sec:subsolY}, $\beta >0$ be small enough so that
$$
\left\{
\begin{array}{l}
0<\beta < \lambda ,\\
\nu (\beta \mbox{Le})-\beta^2 +c\beta \mbox{Le} >0 ,
\end{array}
\right.
$$
and $\gamma > 0$ large enough so that
$$
\left\{
\begin{array}{l}
\displaystyle \gamma \times \min_{\overline{\omega}} \psi_{\beta \text{Le}} \geq 1 ,\\
\displaystyle \gamma \mbox{Le}^{-1} (\nu (\beta \mbox{Le})-\beta^2 +c\beta \mbox{Le}) \times \min_{\overline{\omega}} \psi_{\beta \text{Le}} >  C \max_{y \in \overline{\omega}} \frac{\partial f}{\partial T}(y,0),
\end{array}
\right.
$$
where $\psi_{\beta \text{Le}}$ is the positive eigenfunction of (\ref{eqn:principaleigenvalue}) with $\lambda=\beta \mbox{Le}$, normalized so that $\|\psi_{\beta \text{Le}}\|_{L^\infty (\omega )} =1$. Let $\underline{Y}$ be defined by
$$\underline{Y} (x,y) = \max(0,1-\gamma \psi_{\beta \text{Le}} (y) e^{-\beta (x-ct)}).$$
Note that $\underline{Y} =0$ for $x \leq ct$. Moreover, $\underline{Y}$ satisfies the Neumann boundary conditions on $\partial \Omega$ and when $\underline{Y} >0$, then $x>ct$ and
\begin{eqnarray*}
&& \displaystyle  \underline{Y}_t - \text{Le}^{-1} \Delta \underline{Y} + u(y) \underline{Y}_x + f(y,T)\underline{Y}\\
&& \displaystyle \leq \underline{Y}_t - \text{Le}^{-1} \Delta \underline{Y} + u(y) \underline{Y}_x + f(y,\overline{T})\underline{Y}\\
&& \displaystyle \leq -\gamma \text{Le}^{-1} ( \nu(\beta \text{Le}) -\beta^2 +c\beta \text{Le}) \psi_{\beta \text{Le}} (y) e^{-\beta (x-ct)} + \frac{\partial f}{\partial T}(y,0) C \phi_{\lambda} (y) e^{-\lambda (x-ct)} \\
&& \displaystyle \leq  -\gamma \text{Le}^{-1} ( \nu(\beta \text{Le}) -\beta^2 +c\beta \text{Le}) \psi_{\beta \text{Le}} (y) e^{-\beta (x-ct)} +\frac{\partial f}{\partial T}(y,0) C e^{-\beta (x-ct)} \leq 0 .
\end{eqnarray*}
The maximum principle then implies that $$Y(t,x,y) \geq \underline{Y} (t,x,y)$$ for all $t\geq 0$ and $(x,y)\in \overline{\Omega}$, provided that
$$Y_0 (x,y) \geq 1-\gamma \psi_{\beta \text{Le}} (y) e^{-\beta x},$$
which indeed holds for $\beta \leq \lambda '$ and $\gamma \times \min_{\overline{\omega}} \psi_{\beta \text{Le}} \geq C_3$ (this is possible since $\beta$ could be chosen arbitrarily small and $\gamma$ arbitrarily large).

Lastly, as in Section \ref{sec:subsolT}, the fact that $\lambda$ is the smallest positive root of $k(s)=cs$ allows us to choose $\eta >0$ small enough so that
$$
\left\{
\begin{array}{l}
0<\eta < \min (\beta ,\alpha \lambda) ,\\
\varepsilon :=c(\lambda +\eta )-k(\lambda + \eta )>0 ,
\end{array}
\right.
$$
where $\alpha >0$ such that $f(y,.)$ and $h(y,.)$ are of class $C^{1,\alpha } ([0,s_0 ])$ for some $s_0 >0$ uniformly in $y \in \overline{\omega}$. Let $M \geq 0$ such that
$$
\left\{
\begin{array}{l}
\displaystyle f(y,s)\geq \frac{\partial f}{\partial T} (y,0)s - Ms^{1+\alpha } ,\\
\displaystyle h(y,s)\leq \frac{\partial h}{\partial T}(y,0)s + Ms^{1+\alpha } ,\\
\end{array}
\right. \mbox{ for all } s \in [0,s_0]  \mbox{ and for all } y \in \overline{\omega}.
$$
Now take $\xi \geq 0$ sufficiently large so that
$$\underline{Y} (t,x,y) = 1 -\gamma \psi_{\beta \text{Le}} (y) e^{-\beta (x-ct)} \mbox{ when } x-ct \geq \xi .$$
Next, let $\delta >0$ large enough so that
\begin{equation}\label{eqn:delta}
\left\{
\begin{array}{l}
\displaystyle \phi_{\lambda} (y) e^{-\lambda (x-ct)} - \delta \phi_{\lambda + \eta} (y) e^{-(\lambda +\eta ) (x-ct)} \leq s_0 \mbox{ for all } t\geq 0 \mbox{ and } (x,y)\in\overline{\Omega} ,\\
\displaystyle \phi_{\lambda} (y) e^{-\lambda (x-ct)} - \delta \phi_{\lambda + \eta} (y) e^{-(\lambda +\eta ) (x-ct)} \leq 0 \mbox{ \ when } x-ct \leq \xi  ,\\
\displaystyle \delta \varepsilon \times \min_{\overline{\omega}} \phi_{\lambda_c + \eta } \geq \gamma \max_{y \in \overline{\omega}} \frac{\partial f}{\partial T} (y,0) + 2M .
\end{array}
\right.
\end{equation}
We then define
$$\underline{T} (t,x,y) = C' \max \left( 0 ,  \phi_{\lambda } (y) e^{-\lambda (x-ct)} -\delta \phi_{\lambda + \eta} (y) e^{-(\lambda +\eta ) (x-ct)} \right) ,$$
where $C' \in (0,1)$ to be determined. Note first that $0 \leq \underline{T} \leq s_0$ in $\overline{\Omega}$, and $\underline{T} =0$ when $x-ct \leq \xi$. Thus, if $\underline{T} (x,y) > 0$, then $x-ct > \xi \geq 0$ and $0\leq \underline{Y}(t,x,y)=1-\gamma \psi_{\beta \text{Le}} (y) e^{-\beta (x-ct)}$.

Then, in that case, we have:
\begin{eqnarray*}
&& \displaystyle \underline{T}_t - \Delta \underline{T}  + u(y)\underline{T}_x - f(y,\underline{T})Y + h(y,\underline{T})\\
&& \displaystyle \leq \underline{T}_t - \Delta \underline{T}  + u(y)\underline{T}_x - f(y,\underline{T})\underline{Y} + h(y,\underline{T})\\
&& \displaystyle \leq \underline{T}_t -\Delta \underline{T}  + u(y)\underline{T}_x - (\frac{\partial f}{\partial T} (y,0) \underline{T}-M\underline{T}^{1+\alpha})(1-\gamma \psi_{\beta \text{Le}} (y) e^{-\beta (x-ct)})+ \frac{\partial h}{\partial T} (y, 0)\underline{T}+M\underline{T}^{1+\alpha}\\
&& \displaystyle \leq \delta C' ( k(\lambda +\eta ) -c(\lambda +\eta))\phi_{\lambda +\eta} (y) e^{-(\lambda +\eta) (x-ct)} + \frac{\partial f}{\partial T} (y,0) \gamma \underline{T}\psi_{\beta \text{Le}} (y) e^{-\beta (x-ct)} + 2M \underline{T}^{1+\alpha}\\
&& \displaystyle \leq - \delta \varepsilon C' \phi_{\lambda +\eta} (y) e^{-(\lambda +\eta)(x-ct)} + \frac{\partial f}{\partial T} (y,0)C' \gamma e^{-(\lambda + \beta )(x-ct)} + 2MC' e^{-\lambda (1+\alpha )(x-ct)}\\
&& \displaystyle \leq C' (- \delta \varepsilon \phi_{\lambda +\eta} (y) + \frac{\partial f}{\partial T} (y,0)\gamma + 2M)e^{-(\lambda + \eta )(x-ct)} \leq 0 .
\end{eqnarray*}
Moreover, the function $\underline{T}$ satisfies the Neumann boundary conditions on $\partial \Omega$. In order to apply the maximum principle, it now remains to check that $T_0 (x,y) \geq \underline{T} (0,x,y)$ for all $(x,y) \in \overline{\Omega}$. Indeed, for $x \leq \xi$, then $\underline{T} (0,x,y) = 0 \leq T_0 (x,y)$. On the other hand, for $x > \xi \geq 0$, we have that $T_0 (x,y) \geq C_1 e^{-\lambda x} \geq C' \phi_\lambda (y) e^{-\lambda x} \geq \underline{T} (0,x,y)$, provided that $C' < C_1$.
Therefore, it follows from the maximum principle that
$$T(t,x,y) \geq \underline{T} (t,x,y)$$
for all $t\geq 0$ and $(x,y) \in \overline{\Omega}$. Let now $x_0 \in \mathbb{R}$ such that $$ \alpha (x_0 ,y) := \phi_{\lambda } (y) e^{-\lambda x_0} -\delta \phi_{\lambda + \eta} (y) e^{-(\lambda +\eta ) x_0} > 0$$ for all $y \in \overline{\omega}$. We then have that $T(t,x_0 +ct ,y) \geq  \underline{T} (t,x_0 + ct,y) = \alpha (x_0 ,y) >0$, which concludes the proof of part (c) of Theorem \ref{th:cauchy}. $\Box$

\nocite{hamel-quenching,abel1,berestycki0,berestycki1,berestycki2,berestycki3,berestycki4}
\nocite{constantin1,constantin2,heinze1,khoudier1,kagan1,kagan2,kiselev1,kiselev2,nolen1}
\nocite{vladimirova1,volpert1,xin1,xin2,xin3,berestycki5,ducrot1,collet1}
\nocite{hamel-adiabatic,hamel-nonadiabatic,marion1,berestycki6}

\bibliographystyle{plain}
\bibliography{KPP5}

\begin{thebibliography}{10}

\bibitem{abel1}
M.~Abel, A.~Celani, D.~Vergni, and A.~Vulpiani.
\newblock Front propagation in laminar flows.
\newblock {\em Physical review E}, 64, 2001.

\bibitem{berestycki0}
B.~Audoly, H.~Berestycki, and Y.~Pomeau.
\newblock R{\'e}action-diffusion en {\'e}coulement stationnaire rapide.
\newblock {\em C. R. Acad. Sci. Paris}, 328 II:255--262, 2000.

\bibitem{berestycki1}
H.~Berestycki.
\newblock The influence of advection on the propagation of fronts in
  reaction-diffusion equations.
\newblock {\em Nonlinear PDEs in Condensed Matter and Reactive Flows, NATO
  Science Series C}, 589, 2002.

\bibitem{berestycki2}
H.~Berestycki and F.~Hamel.
\newblock Front propagation in periodic excitable media.
\newblock {\em Comm. Pure Appli. Math}, 55:949--1032, 2002.

\bibitem{berestycki-hamel}
H.~Berestycki and F.~Hamel.
\newblock {\em Reaction-diffusion equations and propagation phenomena}.
\newblock Springer-Verlag, to appear.

\bibitem{hamel-quenching}
H.~Berestycki, F.~Hamel, A.~Kiselev, and L.~Ryzhik.
\newblock Quenching and propagation in {KPP} reaction-diffusion equations with
  a heat loss.
\newblock {\em Arch. Ration. Mech. Anal.}, 178:57--80, 2005.

\bibitem{berestycki5}
H.~Berestycki, F.~Hamel, and N.~Nadirashvili.
\newblock The speed of propagation for {KPP} type problems. {I} - {P}eriodic
  framework.
\newblock {\em J. European Math. Soc.}, 7:173--213, 2005.

\bibitem{berestycki6}
H.~Berestycki, F.~Hamel, and N.~Nadirashvili.
\newblock The speed of propagation for {KPP} type problems. {II} - {G}eneral
  domains.
\newblock {\em J. Amer. Math. Soc.}, to appear.

\bibitem{berestycki3}
H.~Berestycki, B.~Larrouturou, and P.-L. Lions.
\newblock Multi-dimensional traveling wave solutions of a flame propagation
  model.
\newblock {\em Arch. Rational Mech. Anal.}, 111:33--49, 1990.

\bibitem{berestycki8}
H.~Berestycki, B.~Larrouturou, P.-L. Lions, and J.-M. Roquejoffre.
\newblock An elliptic system modelling the propagation of a multidimensional
  flame.
\newblock {\em Unpublished manuscript}, 1995.

\bibitem{berestycki4}
H.~Berestycki and L.~Nirenberg.
\newblock Traveling wave in cylinders.
\newblock {\em Annales de l'IHP, Analyse non lin{\'e}aire}, 9:497--572, 1992.

\bibitem{collet1}
P.~Collet and J.~Xin.
\newblock Global existence and large time asymptotic bounds of ${L}^\infty$
  solutions of thermal diffusive combustion systems on $\mathbb{R}^n$.
\newblock {\em Ann. Scuola Norm. Sup. Pisa Cl. Sci.}, 23:625--642, 1996.

\bibitem{constantin1}
P.~Constantin, A.~Kiselev, A.~Oberman, and L.~Ryzhik.
\newblock Bulk burning rate in passive-reactive diffusion.
\newblock {\em Arch. Rational Mech. Anal.}, 154:53--91, 2000.

\bibitem{constantin2}
P.~Constantin, A.~Kiselev, and L.~Ryzhik.
\newblock Quenching of flames by fluid advection.
\newblock {\em Comm. Pure Appl. Math}, 54:1320--1342, 2001.

\bibitem{ducrot1}
A.~Ducrot.
\newblock Multi-dimensional combustion waves for {L}ewis number close to one.
\newblock {\em Math. Methods Appl. Sci.}, 30:291--304, 2007.

\bibitem{giletti2}
T.~Giletti.
\newblock {KPP} reaction-diffusion system with loss inside a cylinder:
  convergence toward the problem with robin boundary conditions, preprint.

\bibitem{giovangigli1}
V.~Giovangigli.
\newblock Nonadiabatic plane laminar flames and their singular limits.
\newblock {\em SIAM J. Math. Anal.}, 21:1305--1325, 1990.

\bibitem{gordon-periodic}
P.~Gordon, L.~Ryzhik, and N.~Vladimirova.
\newblock The {KPP} system in a periodic flow with a heat loss.
\newblock {\em Nonlinearity}, 18:571--589, 2005.

\bibitem{hamel-nonadiabatic}
F.~Hamel and L.~Ryzhik.
\newblock Non-adiabatic {KPP} fronts with an arbitrary lewis number.
\newblock {\em Nonlinearity}, 18:2881--2902, 2005.

\bibitem{hamel-adiabatic}
F.~Hamel and L.~Ryzhik.
\newblock Travelling waves for the thermodiffusive system with arbitrary lewis
  numbers.
\newblock {\em Arch. Ration. Mech. Anal.}, preprint.

\bibitem{hamel1}
F.~Hamel and Y.~Sire.
\newblock Spreading speeds for some reaction-diffusion equations with general
  initial conditions, preprint.

\bibitem{heinze1}
S.~Heinze, G.~Papanicolaou, and A.~Stevens.
\newblock Variational principles for propagation speeds in inhomogeneous media.
\newblock {\em SIAM J. Appl. Math}, 62:129--148, 2001.

\bibitem{kagan2}
L.~Kagan, P.D. Ronney, and G.~Sivashinsky.
\newblock Activation energy effect on flame propagation in large-scale vortical
  flows.
\newblock {\em Combust. Theory Modelling}, 6:479--485, 2002.

\bibitem{kagan1}
L.~Kagan and G.~Sivashinsky.
\newblock Flame propagation and extinction in large-scale vortical flows.
\newblock {\em Combust. Flame}, 120:222--232, 2000.

\bibitem{khoudier1}
B.~Khoudier, A.~Bourlioux, and A.~Majda.
\newblock Parametrizing the burning rate speed enhancement by small scale
  periodic flows: I. {U}nsteady shears, flame residence time and bending.
\newblock {\em Combust. Theory Model.}, 5:295--318, 2001.

\bibitem{kiselev1}
A.~Kiselev and L.~Ryzhik.
\newblock Enhancement of the travelling front speeds in reaction-diffusion
  equations with advection.
\newblock {\em Ann. Inst. H. Poincar{\'e} Anal. Non Lin{\'e}aire}, 18:309--358,
  2001.

\bibitem{kiselev2}
A.~Kiselev and L.~Ryzhik.
\newblock An upper bound for the bulk burning rate for systems.
\newblock {\em Non linearity}, 14:1297--1310, 2001.

\bibitem{marion1}
M.~Marion.
\newblock Qualitative properties of a nonlinear system for laminar flames
  without ignition temperature.
\newblock {\em Nonlinear Anal. Th. Meth. Appl.}, 9:1269--1292, 1985.

\bibitem{murray1}
J.D. Murray.
\newblock {\em Mathematical biology. {I} {A}n introduction}.
\newblock Springer, third edition, 2002.

\bibitem{murray2}
J.D. Murray.
\newblock {\em Mathematical biology. {II} {S}patial models and biomedical
  applications}.
\newblock Springer, third edition, 2003.

\bibitem{nolen1}
J.~Nolen and J.~Xin.
\newblock Reaction-diffusion front speeds in spatially-temporally periodic
  shear flows.
\newblock {\em SIAM Jour. MMS}, 1:554--570, 2003.

\bibitem{roques1}
L.~Roques.
\newblock Existence de deux solutions du type front progressif pour un
  mod{\`e}le de combustion avec pertes de chaleur.
\newblock {\em C. R. Acad. Sci. Paris Ser. I}, 340:493--497, 2005.

\bibitem{roques2}
L.~Roques.
\newblock Study of the premixed flame model with heat losses the existence of
  two solutions.
\newblock {\em Euro. J. Appl. Math.}, 16:741--765, 2005.

\bibitem{vladimirova1}
N.~Vladimirova, P.~Constantin, A.~Kiselev, O.~Ruchayskiy, and L.~Ryzhik.
\newblock Flame enhancement and quenching in fluid flows.
\newblock {\em Combust. Theory Model.}, 7:487--508, 2003.

\bibitem{volpert1}
V.A. Volpert and A.I. Volpert.
\newblock Existence and stability of multidimensional travelling waves in the
  monostable case.
\newblock {\em Israel Jour. Math.}, 110:269--292, 1999.

\bibitem{xin1}
J.~Xin.
\newblock Existence of planar flame fronts in convective-diffusive periodic
  media.
\newblock {\em Arch. Rat. Mech. Anal.}, 121:205--233, 1992.

\bibitem{xin2}
J.~Xin.
\newblock Existence and nonexistence of travelling waves and reaction-diffusion
  front propagation in periodic media.
\newblock {\em Jour. Stat. Phys.}, 73:893--926, 1993.

\bibitem{xin3}
J.~Xin.
\newblock Analysis and modelling of front propagation in heterogeneous media.
\newblock {\em SIAM Rev.}, 42:161--230, 2000.

\end{thebibliography}

\end{document}